\documentclass[a4paper,10pt]{article}
\usepackage[english]{babel}

\usepackage[a4paper,left=3cm,right=3cm,top=2.5cm,bottom=2.5cm]{geometry}

\usepackage{hyperref}
\usepackage{bm}

\usepackage{authblk}

\usepackage{siunitx}
\usepackage{graphicx}
\usepackage{tikz}
\usepackage{transparent}
\usepackage{amsfonts}
\usepackage{mathtools}
\usepackage{amsmath}
\usepackage{subcaption}
\graphicspath{ {./images/} }

\title{A two stages Deep Learning Architecture for Model Reduction of Parametric Time-Dependent Problems}
\author[1]{Isabella Carla Gonnella \thanks{igonnell@sissa.it}}
\author[1]{Martin W. Hess \thanks{mhess@sissa.it}}
\author[1,2]{Giovanni Stabile \thanks{giovanni.stabile@uniurb.it (Present Affiliation: University of Urbino)}}
\author[1]{Gianluigi Rozza \thanks{grozza@sissa.it}}

\affil[1]{Mathematics Area, mathLab, SISSA, via Bonomea 265, I-34136 Trieste, Italy}
\affil[2]{Department of Pure and Applied Sciences, Informatics and Mathematics Section, University of Urbino Carlo Bo, Piazza della Repubblica, 13, I-61029 Urbino, Italy}
\date{} 

\begin{document}

\maketitle

\begin{abstract}
Parametric time-dependent systems are of a crucial importance in modeling real phenomena, often characterized by non-linear behaviours too. Those solutions are typically difficult to generalize in a sufficiently wide parameter space while counting on limited computational resources available. As such, we present a general two-stages deep learning framework able to perform that generalization with low computational effort in time. It consists in a separated training of two pipe-lined predictive models. At first, a certain number of independent neural networks are trained with data-sets taken from different subsets of the parameter space. Successively, a second predictive model is specialized to properly combine the first-stage guesses and compute the right predictions. Promising results are obtained applying the framework to incompressible Navier-Stokes equations in a cavity (Rayleigh-Bernard cavity), obtaining a $97\%$ reduction in the computational time comparing with its numerical resolution for a new value of the Grashof number.
\end{abstract}

\providecommand{\keywords}[1]
{
  \small    
  \textbf{\textit{Keywords---}} #1
}
\keywords{
reduced order modeling, deep learning , long-short term memory networks , convolutional layers , time forecasting , time-dependent parametric PDEs}

\date{} 

\section{Introduction}
Time-dependent systems, especially in the parametrized setting, describe a huge number of problems and are therefore a pervasive topic of extended scientific interest and industrial value. Indeed, \emph{parametric dynamical systems} modeling and control play a fundamental role in many research fields, as in the case of fluid dynamics, chemical reactions, biological problems and more. 

In the majority of scenarios, the most suitable way to study such dynamics passes through numerical simulation. Especially for what concerns problems modelled by differential and partial differential equations, numerical approximation represent the standard to compute the system's response. 

However, a problem of dimensionality of the system's numerical discretization often appears significant, as performing multiple simulations in large-scale settings typically reveals demands of computational resources difficult to handle. 

This gives rise to the need of finding alternatives to classical numerical methods (Finite Element Method, Finite Volume Method, Finite Difference Method) in order to approximate the parametric response of a given system at a reduced computational cost. Reduced order models (ROMs) demonstrated to be a powerful tools in this regard and nowadays it is possible to find a large variety of applications in a number of different fields as heat transfer, fluid dynamics, shape optimization, uncertainty quantification. The main idea of ROMs is to approximate a high dimensional model, usually referred as full order model (FOM), with a low dimensional one still preserving the solution's key features. There mainly exist two different techniques to obtain a ROM: intrusive and non-intrusive approaches. The common feature of both approaches is the computational splitting into two distinct phases: an offline (or training), where the parametric response of the system is explored for selected values of the input parameters, and online (or testing) one that allow to retrieve the system's response for any new value of the input parameters \cite{aroma_book}. In both cases the results acquired during the initial exploration of the solution manifold are used to perform a compression of the discrete solution manifold. It can be performed using both linear (proper orthogonal decomposition, reduced basis methods) or nonlinear approaches (autoencoders, convolutional autoenconders). The two differ in the methodology used to approximate the evolution of the latent coordinates (reduced basis coefficients) in the latent space (reduced basis space). 

Intrusive methods, that have its root in the classical field of scientific computing, use a Galerkin-(Petrov) projection of the system of equations describing the dynamics onto a linear-(nonlinear) subspace-(manifold) in order to generate a low dimensional model that need to be solved for any new value of the input parameters. These techniques, exploiting the underlying physical principles generally exhibit better generalization properties and perform well with less training data (\cite{HijaziStabileMolaRozza2020b,GeorgakaStabileRozzaBluck2019,RomorStabileRozza2022}). On the other hand, they show severe limitations when addressing nonlinear time-dependent parametric PDEs, due in general to the difficulty of capturing complex physical patters and generalizing them to a large set of \emph{online parameters} \cite{Quarteroni2015ReducedBM,aroma_book}.

Non-intrusive approaches are instead solely based on input-output data and do not require the explicit knowledge of the underlying equations. The evolution of the latent coordinates is retrieved by means of different regression or interpolation techniques. Being \emph{data-driven} they have the significant advantage of making the methods \textit{non-intrusive}, \textit{i.e.} allowing the high-fidelity model to be run in “black-box” mode, needing only to specify a set of input parameters and generate the corresponding system outputs. In this article we will focus only on the second type of methods (i.e. non-intrusive methods) and particularly on approaches suitable to address parameter and time dependent problems. Many research works that employ non-intrusive methods are in fact dedicated to stationary parameter-dependent problems or to transient problems \cite{CAE_CFD}, but far less material is available for transient and parameter-dependent problems.

About parameter-dependent problems, some developments are present dealing with properly enhanced reduced basis methods \cite{Hesthaven_2022}. In particular, the use of data-driven techniques shows itself to be a key tool in the formulation of reduced basis methods that are both stable and highly efficient, even for general nonlinear problems. This is achieved by introducing non-intrusive reduced order models in which a data-driven map is learned as the map between parameter space and coefficients of the reduced basis to reconstruct the solutions, as in the case of \cite{PIN_ROM,non_intrusive_NN,ROM_gaussianp}. However, these approaches do not deal with time-dependent problems.

An example of an approach in that direction is provided by \cite{HessQuainiRozza_ACOM_2022}, where a combination of Proper Orthogonal Decomposition (POD), Dynamic Mode Decomposition (DMD) and Manifold Interpolation is developed to approximate a given time-trajectory. An other case in which POD is utilized and also time dependence is considered, is found in \cite{ROM_NN_time}, where time is treated as an extra parameter.

Instead, for what concerns machine learning techniques, many of them have demonstrated to be particularly useful in the approximation of nonlinear dynamics. It is the case of models such as SVM \cite{SVM}, ARIMA \cite{ARIMA}, as well as probabilistic ones involving hidden Markov models \cite{HM} or fuzzy logic \cite{fuzzy}. Finally, Artificial Neural Networks (ANNs) have been recently massively considered to provide fast and reliable approximations of PDE solutions, thanks to the universal approximation theorem \cite{Hornik1991} that led to different proposals on the topic \cite{ODE_NN}\cite{NeuralDE}\cite{seq2seq}. Specifically, ANNs provided with internal recurrence mechanisms have gradually become the standard for time series prediction when dealing with large amount of data available for the training \cite{time_window_LSTM}\cite{LSTM_chaos}\cite{LSTMvsARIMA}.

However, ANNs express interesting potentialities not only for what concerns sequential learning with memory-aware networks, but also with tools to operate nonlinear dimensionality reduction such as Convolutional Auto-Encoders (CAE), which are actually employed in many recent works \cite{ROM_CAE, DL-ROM,Muecke2021}. These works actually deal with time-dependent systems, thus including some kind of time prediction methodologies after the first nonlinear reduction with CAEs. For instance, in the first reference the time stepping is done intrusively using multistep methods on the reduced model derived from a Galerkin projection procedure, while in the last ones LSTMs and FFNNs are used for time stepping of the reduced state.  Multi-level CAEs are moreover used in \cite{Multi_levelCAE}, employed to reduce the spatial and temporal dimensions of the problem. In addition, POD and CAE are sometimes used one after the other in the same dimensionality reduction process, as in the case of \cite{Cracco2022}.

It is to be noted that much of the success of Artificial Neural Networks (ANN) based ROM has been boosted further by the availability of open source software frameworks such as PyTorch \cite{pytorch} and Tensorflow \cite{tensorflow}. Indeed, they have made implementation and training possible without expert knowledge, also exploiting the eventual availability of computation accelerating hardware such as GPUs, which has made training of very large models feasible.

In this work a novel two-stages memory-aware ANNs model order reduction approach is developed. At the best of our knowledge, it implements a different strategy with respect to what has been already proposed in the field of trainable architectures able to generalize parametric time-dependent dynamics with scarce set of available solutions.

A windowed approach involving LSTMs (see Appendix \ref{appendix}) for the time-stepping is chosen, meaning that, given a time series forecasting problem, we aim to find: $$\mathcal{J}: (f(x_{t-p};\mathbf{\theta}),\dots,f(x_t;\mathbf{\theta})) \longrightarrow (f(x_{t+1};\mathbf{\theta}),\dots,f(x_{t+m};\mathbf{\theta}))$$ where the time-series $f(\cdot)$ is dependent on the parameters $\mathbf{\theta}$, being $ (f(x_{t-p};\mathbf{\theta}),\dots,f(x_t;\mathbf{\theta}))$ the input time-window. Windowed regressive networks \cite{time_window_LSTM} have been already exploited in multiple applications such as \emph{neural ODEs} \cite{neural_ODE}, where deterministic numerical solvers are led to consider also statistically learned residuals to perform the PDE integration, but also in \cite{LSTM_sea}, where a time series approach using LSTMs proved to be effective in forecasting the sea surface temperature in marine systems. It is in general to be noted that an architecture aimed to find an effective correlation between past sequenced and future ones exhibits close similarities with the behavior of numerical solvers: both of them build predictions for future times based on a certain number of the past ones.

More in depth, differently from what has been done until now, our Neural Network architecture implements a \emph{partitioning-averaging} approach to the parametric problem. It requires different models to be trained for different areas of the parameter space. Their predictions are subsequently combined in a weighted proper way depending on the new parameter for which the prediction is asked. This strategy has in principle the advantage to be able to learn an internal non-linear representation of the qualitative changes of a system with respect to the action of a certain set of parameters. Indeed, it breaks in two parts the reproduction of multiple potentially different local dynamics and their generalization to any new parameter belonging to the considered space.

LSTM-derived neural networks are used in a two stages framework (described in Section \ref{method}) for their ability in learning both short and long-time dependencies in the data, which make them particularly important among all the different recurrent cells (see Appendix \ref{appendix}). Moreover, with this architecture an arbitrary long prediction in time can be obtained thanks to an auto-sustained iterative mechanism, that updates each time the input of the framework with the previous predictions.

Such framework generalization capabilities have been firstly tested on ODEs systems, whose results are available in Section \ref{section:ODE}, and secondly on a widely used benchmark considering the incompressible Navier-Stokes equations in a rectangular cavity (see Section \ref{section:cavity}): the Rayleigh-Benard cavity problem. In order to deal with such a high-dimensional discredized system, the example reported in \cite{Fresca2022} with the POD-DL-ROM has been followed, and a POD has been previously performed to reduce the dimensions, speeding up the training phase.

This last test case considers as the model parameter the Grashof number $Gr$, which is a non-dimensional quantity that describes the ratio of buoyancy forces to viscous forces. It is to be noted that, although this problem considers only one physical parameter,it exhibits a wide range of patterns. Indeed, if at low Grashof numbers the system has unique steady-state solutions, as $Gr$ increases the system undergoes several Hopf bifurcations and multiple solutions arise for the same parameter value. Such solutions past the Hopf bifurcations result time-dependent, being time-periodic at medium Grashof numbers, and exhibiting turbulent
behaviour at very high $Gr$ values. A particular difficulty in applying a ROM approach to the Rayleigh-Benard cavity over a large range of Grashof numbers is related to the fact that frequencies of time-periodic solutions could significantly vary in such range, making hard an exact approximation of the solution for a general \emph{online parameter}.

Our tests apply the new model reduction approach to a range of $50\cdot 10^3$ medium Grashof numbers, taking as parameter space the interval $Gr \in \{100\cdot 10^3,150\cdot 10^3\}$.

\newpage

\section{Methodology}\label{method}

\subsection{Two-stages architecture}\label{architecture}
The proposed \textit{data-driven} approach is realized through a two-stages architecture, which can be interpreted as an implementation of a \emph{partitioning-averaging} method, trained to potentially reveal the system's non-linear dependencies on the considered parameters. The partitioning-averaging method generates accurate estimations valid over local partitions in the first stage, while the second one globally averages the local estimates in an appropriate sense.
This approach implements a regression method based on \emph{k-means} clustering \cite{k_means}, which is a standard method to cluster data vectors.

Here, the \emph{k-means} clustering is performed in the sampled parameter space $\Theta_{\textrm{training}} = \{\mathbf{\theta}_{t}^i\}_{i=1}^n$, which is assumed to be a sufficiently fine sample of the $p$-dimensional parameter space $\Theta \subset \mathbb{R}^p$. The k-means clustering results in $k$ different data sets (or clusters), which form a partition of $\Theta_{\textrm{training}}$. The \emph{centroid} of each cluster is denoted $\{\mathbf{\theta}_{c}^i\}_{i=1}^k$.
Each parameter vector $\mathbf{\theta}_{t}^i \in \Theta_{\textrm{training}}$ defines a trajectory
\begin{align} 
 \{ \mathbf{x}_{1}^{i}, \mathbf{x}_{2}^{i},\dots, \mathbf{x}_{T}^{i} \}\ \ \ \ \ \textrm{where}\ \ \ \ \ \mathbf{x}_j^i \in \mathbb{R}^z\ \ \ \ \ \forall j \in \{ 1,\dots,T \}
\end{align}
through the solution of the respective ODE or PDE, where $z$ represents the number of variables evolving and $T$ the number of time steps. The solution trajectories corresponding to the training parameter values of the same cluster are concatenated, forming the final data-sets $\{D_i\}_{i=1}^k$.

\begin{figure}[b!]
    \centering
    \input{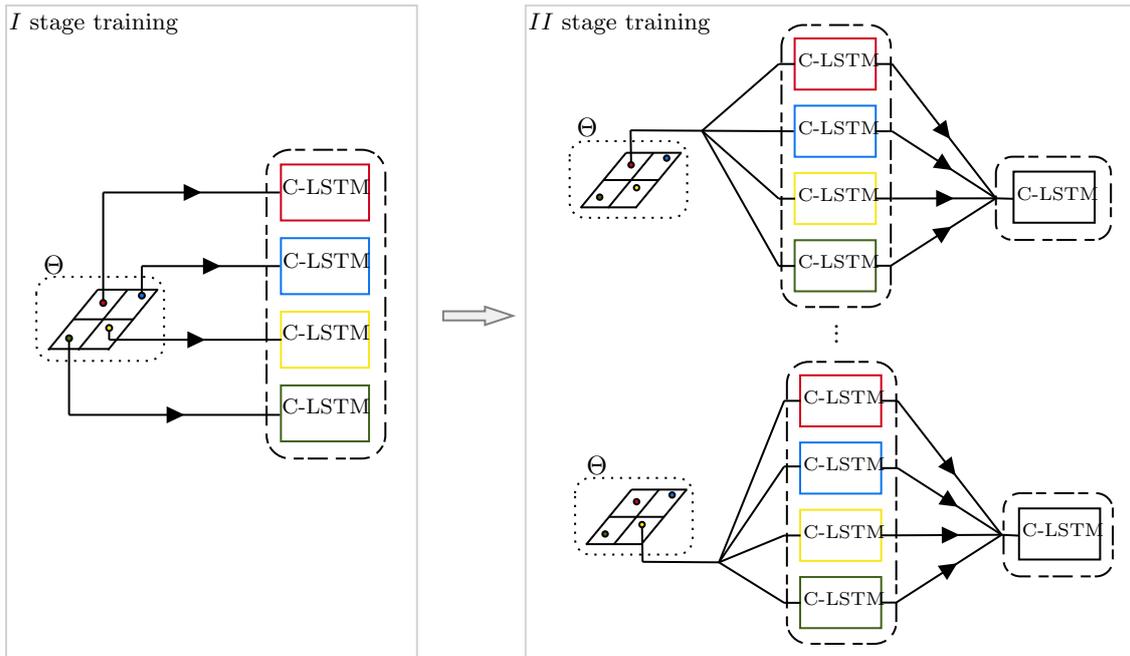}
    \caption{\small Example of the training phases of the two stages assuming a 2-dimensional parameter space $\Theta$ with $n=k=4$.}
    \label{fig:training}
\end{figure}

Subsequently, the $k$ Neural Networks (NNs) of the first stage are trained respectively with those $k$ generated data-sets $\{D_i\}_{i=1}^k$, resulting in a set of $k$ \emph{localized} models (see Figure \ref{fig:training}). More precisely, we can approximate the trained models with a set of $k$ functions: $$\{\mathcal{F}_i(\mathbf{x}_{t-w+1},\mathbf{x}_{t-w+2}, \dots,\mathbf{x}_{t};\mathbf{\theta)}\}_{i=1}^k \ \ \ \ \ \textrm{s.t.}\ \ \ \ \  \mathcal{F}_i: (\mathbb{R}^{w} \times\ \mathbb{R}^z; \mathbb{R}^p)\rightarrow \mathbb{R}^m \textrm{x}\ \mathbb{R}^z,$$ where $w$ is the size of the past system evolution time-window, while $m$ represents the number of next time-step predictions about the system dynamics that the model has been trained to perform.

Coming to the second stage of the architecture, here a Neural Network receives as input all the outputs $\{f_i\}_{i=1}^k$ of the $k$ first-stage models, and aims to implement an ``averaging-function’’ between these first ``local’’ predictions, based on the difference between the respective \emph{centroids} $\{\mathbf{\theta}_{c}^i\}_{i=1}^k$ and the current parameter values $\mathbf{\theta}$:  $$\mathcal{G}(f_1,f_2,\dots,f_k,\mathbf{\theta}_{c}^1,\mathbf{\theta}_{c}^2,\dots,\mathbf{\theta}_{c}^k; \mathbf{\theta}).$$

\begin{figure}[h!]
    \centering
    \tikzset{every picture/.style={line width=0.75pt}} 

\begin{tikzpicture}[x=0.75pt,y=0.75pt,yscale=-1,xscale=1]

\draw   (38,37.24) -- (177.44,37.24) -- (177.44,53.26) -- (38,53.26) -- cycle ;
\draw   (237.69,19) -- (323.76,19) -- (323.76,71.5) -- (237.69,71.5) -- cycle ;
\draw  [color={rgb, 255:red, 0; green, 0; blue, 0 }  ,draw opacity=1 ][fill={rgb, 255:red, 74; green, 144; blue, 226 }  ,fill opacity=0.3 ] (384.01,38.62) -- (459.87,38.62) -- (459.87,54.37) -- (384.01,54.37) -- cycle ;
\draw   (30,160.74) -- (178.3,160.74) -- (178.3,176.76) -- (30,176.76) -- cycle ;
\draw   (238.55,142.5) -- (324.62,142.5) -- (324.62,195) -- (238.55,195) -- cycle ;
\draw    (82.5,37.51) -- (82.5,53.26) ;
\draw    (134.92,37.51) -- (134.92,53.26) ;
\draw    (86.73,161.02) -- (86.73,176.76) ;
\draw    (177.44,44.7) -- (234.83,44.7) ;
\draw [shift={(236.83,44.7)}, rotate = 180] [color={rgb, 255:red, 0; green, 0; blue, 0 }  ][line width=0.75]    (10.93,-3.29) .. controls (6.95,-1.4) and (3.31,-0.3) .. (0,0) .. controls (3.31,0.3) and (6.95,1.4) .. (10.93,3.29)   ;
\draw    (178.3,169.03) -- (235.69,169.03) ;
\draw [shift={(237.69,169.03)}, rotate = 180] [color={rgb, 255:red, 0; green, 0; blue, 0 }  ][line width=0.75]    (10.93,-3.29) .. controls (6.95,-1.4) and (3.31,-0.3) .. (0,0) .. controls (3.31,0.3) and (6.95,1.4) .. (10.93,3.29)   ;
\draw    (323.76,44.7) -- (381.15,44.7) ;
\draw [shift={(383.15,44.7)}, rotate = 180] [color={rgb, 255:red, 0; green, 0; blue, 0 }  ][line width=0.75]    (10.93,-3.29) .. controls (6.95,-1.4) and (3.31,-0.3) .. (0,0) .. controls (3.31,0.3) and (6.95,1.4) .. (10.93,3.29)   ;
\draw    (324.62,169.03) -- (382.01,169.03) ;
\draw [shift={(384.01,169.03)}, rotate = 180] [color={rgb, 255:red, 0; green, 0; blue, 0 }  ][line width=0.75]    (10.93,-3.29) .. controls (6.95,-1.4) and (3.31,-0.3) .. (0,0) .. controls (3.31,0.3) and (6.95,1.4) .. (10.93,3.29)   ;
\draw  [dash pattern={on 4.5pt off 4.5pt}]  (416.72,54.37) .. controls (422.71,163.51) and (142.05,58.47) .. (135,159.48) ;
\draw [shift={(134.92,161.02)}, rotate = 272.38] [color={rgb, 255:red, 0; green, 0; blue, 0 }  ][line width=0.75]    (10.93,-3.29) .. controls (6.95,-1.4) and (3.31,-0.3) .. (0,0) .. controls (3.31,0.3) and (6.95,1.4) .. (10.93,3.29)   ;
\draw  [fill={rgb, 255:red, 74; green, 144; blue, 226 }  ,fill opacity=0.3 ] (108.58,160.19) -- (178.3,160.19) -- (178.3,176.76) -- (108.58,176.76) -- cycle ;
\draw    (134.92,161.02) -- (134.92,176.76) ;
\draw    (409.49,38.34) -- (409.49,54.09) ;
\draw    (427.56,38.34) -- (427.56,54.09) ;
\draw   (384.87,162.12) -- (476.33,162.12) -- (476.33,177.32) -- (384.87,177.32) -- cycle ;
\draw    (424.95,161.84) -- (424.95,177.59) ;
\draw    (442.17,161.84) -- (442.17,177.59) ;

\draw (97.75,41.71) node [anchor=north west][inner sep=0.75pt]   [align=left] {$\displaystyle \dotsc $};
\draw (98.14,165.94) node [anchor=north west][inner sep=0.75pt]   [align=left] {$\displaystyle \dotsc $};
\draw (40.14,41.61) node [anchor=north west][inner sep=0.75pt]  [font=\tiny] [align=left] {$\displaystyle x_{T-w+1}$};
\draw (150.19,41.89) node [anchor=north west][inner sep=0.75pt]  [font=\tiny] [align=left] {$\displaystyle x_{T}$};
\draw (31.24,165.94) node [anchor=north west][inner sep=0.75pt]  [font=\tiny] [align=left] {$\displaystyle x_{T-w+1+m}$};
\draw (140.47,165.39) node [anchor=north west][inner sep=0.75pt]  [font=\tiny] [align=left] {$\displaystyle x_{T+m}$};
\draw (410.61,46.16) node [anchor=north west][inner sep=0.75pt]   [align=left] [font=\tiny] {$\displaystyle \dotsc $};
\draw (428.26,43.16) node [anchor=north west][inner sep=0.75pt]  [font=\tiny] [align=left] {$\displaystyle x_{T+m}$};
\draw (383.64,43.16) node [anchor=north west][inner sep=0.75pt]  [font=\tiny] [align=left] {$\displaystyle x_{T+1}$};
\draw (425.08,169.2) node [anchor=north west][inner sep=0.75pt]   [align=left] [font=\tiny] {$\displaystyle \dotsc $};
\draw (442.08,166.39) node [anchor=north west][inner sep=0.75pt]  [font=\tiny] [align=left] {$\displaystyle x_{T+2m}$};
\draw (384.94,166.5) node [anchor=north west][inner sep=0.75pt]  [font=\tiny] [align=left] {$\displaystyle x_{T+m+1}$};
\draw (245.51,30.41) node [anchor=north west][inner sep=0.75pt]  [font=\footnotesize] [align=left] {\begin{minipage}[lt]{46.71pt}\setlength\topsep{0pt}
\begin{center}
Two-stages \\C-LSTM
\end{center}

\end{minipage}};
\draw (246.79,153.31) node [anchor=north west][inner sep=0.75pt]  [font=\footnotesize] [align=left] {\begin{minipage}[lt]{46.71pt}\setlength\topsep{0pt}
\begin{center}
Two-stages \\C-LSTM
\end{center}

\end{minipage}};

\end{tikzpicture}
    \caption{\small Given a time instant $T$, the evolution precedes predicting the next $m$ values of the variable interested per cycle. }
    \label{fig:iteration}
\end{figure}
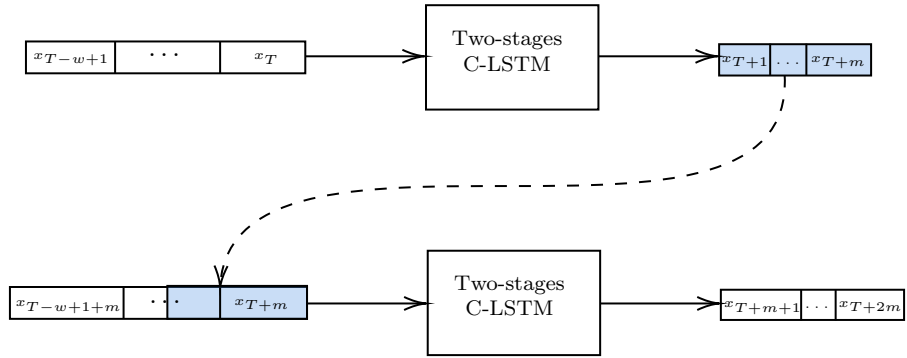

Therefore, the trained architecture can provide an approximation of the time-evolution corresponding to a general \emph{online parameter}, obtained by simply giving as input the first exact time-window $W$ and the new parameter value $\mathbf{\theta_{new}}$. Indeed the evolution is achieved through an iterative recursion, in which the outputs of the architecture are suddenly reused as inputs for the next cycles (Figure \ref{fig:iteration}). 

It is to be noted that, in this way, the advancement in time of the system's variables for a general \emph{online parameter} can be potentially obtained for each desired amount of time-steps, independently from the training solutions' extension. 

Summing up, this framework could be seen as a variant of the Random Forest method \cite{random_forest}, as it builds $k$ different models with $k$ different training data-sets, whose guesses are ``averaged’’ to obtain the final prediction. 

On the other hand, the choice of the data-sets is not ``random’’, but derived from locality-based considerations implemented through the \emph{k-means} algorithm. Hence, a more suitable already proposed methodology, of which our two-stages framework could be considered a generalization, is the \emph{weighted k-Nearest Neighbour} technique \cite{W_KNN}. Indeed, \emph{k-NN} considers the samples $\{\mathbf{\theta}_{t}^i\}_{i=1}^n$ in $\Theta$, and each time a new parameter's value $\mathbf{\theta_n}$ is introduced, its correspondent prediction is computed as the weighted sum of the values associated to its $k$-Nearest Neighbours in the parameters' space.

The differences between what we propose and the \emph{k-NN} method lie in two principal points. Firstly, in our case the values given to the ``weighted averaging function’’ are not the same ones associated to the $k$ $\{\mathbf{\theta}_{c}^i\}_{i=1}^k$, but are computed for the new $\mathbf{\theta_n}$ by the $k$ different models. In second place, to average those values a non-linear function is found by the second-stage NN, being much more complex than a simple weighted average.

Furthermore, it is to be noted that the presented approach is markedly different from other partition-based methods such that the one proposed in \cite{PDE_bif_kmeans}. Indeed, our procedure considers the  k-means employed in the parameter space, not in the space of the discrete PDE solutions.

\subsection{C-LSTM} 
As stated above, both the two stages of the architecture are realized through the exploitation of a particular type of \emph{integrated Long Short Term Memory networks}: the C-LSTM architecture \cite{LSTM-architectures}. 

It consists in a succession of a convolutional and an LSTM layer, two mainstream architectures for such modeling tasks. The usage of the first is aimed at extracting a sequence of higher-level representations of the input, that are successively fed into a recurrent neural network (LSTM) to obtain the final outputs. Indeed, LSTM layers allow to learn from the extracted features' evolution the correct predictions, according with the maintaining of a memory of their long-time and short-time dependencies (appendix \ref{appendix}).

The combination of convolutional neural network (CNN) and LSTM results in a powerful tool for our purposes. In fact, CNN is able to learn local context from temporal or spatial data but lacks the ability of learning sequential correlations. On the other hand, LSTM is specialized for sequential modelling, despite being unable to extract features in a parallel way.

\begin{figure}[h!]
    \centering
    \input{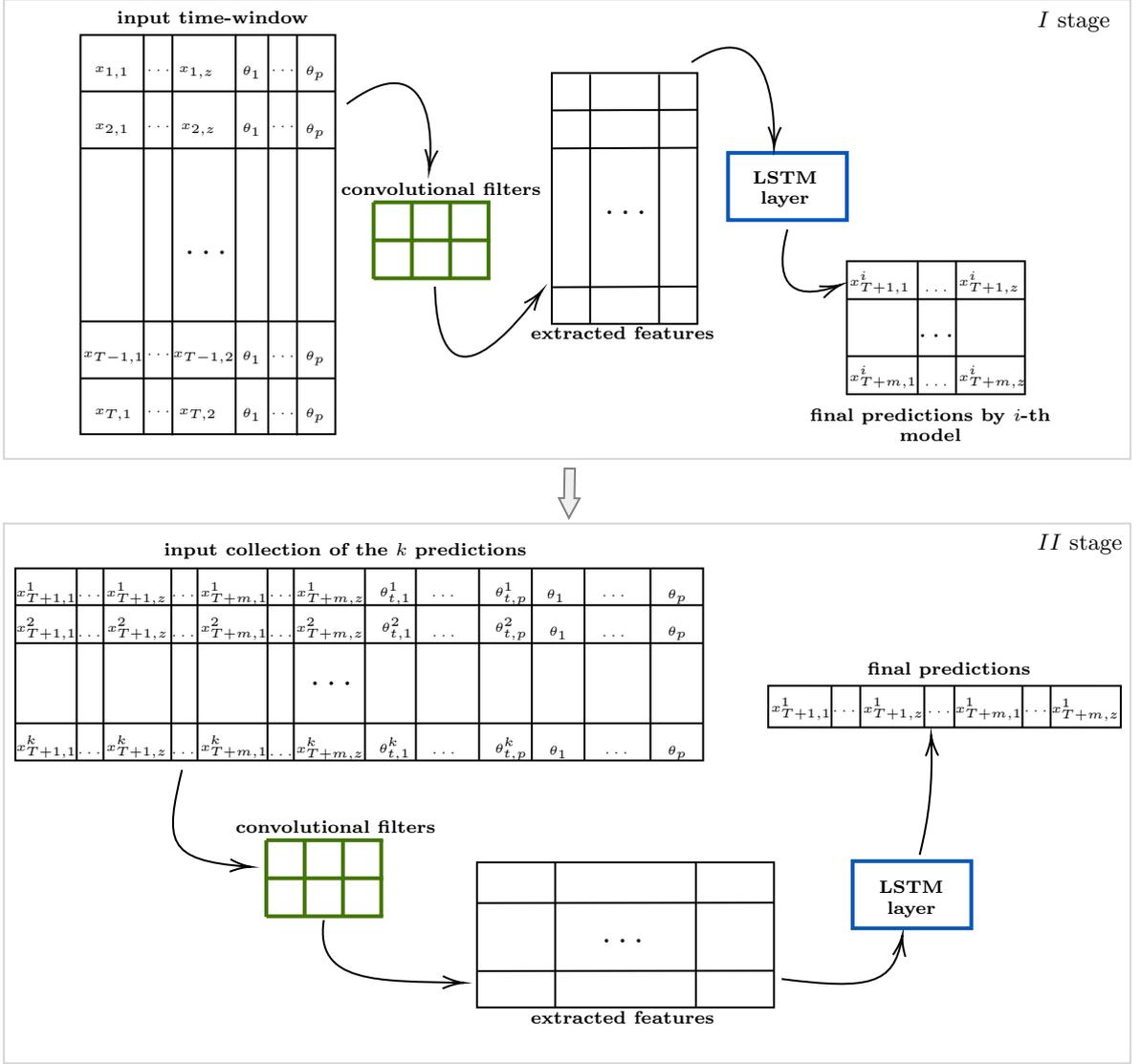}
    \caption{\small General inputs and outputs of the first-stage and second-stage C-LSTM networks, considering general parameters' values $\{\mathbf{\theta}_{i}\}_{i=1}^p$}. 
    \label{fig:stages}
\end{figure}

Examples of C-LSTM employment can be already found in some computer vision tasks like text classification \cite{C-LSTM_text_classification}, image caption \cite{C-LSTM_image_caption} and speech recognition \cite{C-LSTM_speech}.

In our particular case, we need to distinguish between the first stage and the second one, owing to different inputs passed to the two C-LSTM networks (Figure \ref{fig:stages}).

Indeed, while in the first stage the C-LSTM is trained to extract temporal dependencies from the time-window given in input $\{\mathbf{x}_{t-w+1},\mathbf{x}_{t-w+2}, \dots,\mathbf{x}_{t}\}_{i=1}^k$, in the second one it has to learn the spatial dependencies of the $k$ first-stage predictions according to the relation between the respective spatial training parameters $\{\mathbf{\theta}_{c}^i\}_{i=1}^k$ and the current parameter of interest $\mathbf{\theta_{new}}$.

Hence, the combined effect of the two stages results in a pipeline, which is shown in Fig.\ref{fig:stages}. This figure complements Fig.\ref{fig:iteration} in the sense that it provides a detailed view of a single iteration. At first, time-dependencies of the considered variables' trajectories are analysed: a time-window $W$ of the past $w$ temporal steps of such variables is given in input asking to the network to predict their next $m$-steps evolution according to their previous values in $W$ and to the parameters value $\{\mathbf{\theta}_{new_i}\}_{i=1}^p$. In second place, the $k$ outputs of the $k$ first-stage networks are collected and given as input to the second stage. Here, their spatial dependencies are taken into consideration through the extraction of features by the CNN layer (basing on the local distance between the training parameters and the current one), thus elaborated by a LSTM layer.


\section{Application to ODEs} \label{section:ODE}

The proposed architecture has been at first tested on simple ODEs systems in order to prove its generalization capabilities and to investigate the role of some of its parameters, \emph{i.e.} the number of time-steps predicted per iteration and the time-window length.

In particular, we report some results about two examples: the \emph{Duffing Oscillator}, parameterized in its non-linear component and taken with null driving force (\ref{eq:do}), and the \emph{Predator Prey} system in the case of limited resources, with parameterization applied to the predators' growth component (\ref{eq:pp}):
\newline

\begin{minipage}{.45\textwidth}
\begin{equation}
\centering
\left\{\begin{split}
\frac{dr}{dt}=& \ v \\
\frac{dv}{dt}=& \ r -a\cdot r^3 \\ 
\end{split}\right.
\label{eq:do}
 \end{equation}
\end{minipage}
 \begin{minipage}{.45\textwidth}
 \begin{equation}
\centering
\left\{\begin{split}
\frac{dr}{dt}=& \ r\cdot(1-r) - r\cdot v \\
\frac{dv}{dt}=& -v+a\cdot r\cdot v \\ 
\end{split}\right.
\label{eq:pp}
 \end{equation}
\end{minipage}
\newline

The above predictions have been computed considering $w=200$ and $m=1$. The training phase has been performed in the parameter range of $a \in [1,10]$ and $a \in [1,5]$ for the Duffing Oscillator and for the Predatory Prey system respectively, with $k=10$ (training set of $\{\mathbf{\theta}_{c}^i\}_{i=1}^{10} = \{1,2,3,4,5,6,7,8,9,10\}$), and with $k=5$ (training set of $\{\mathbf{\theta}_{c}^i\}_{i=1}^5 = \{1,2,3,4,5\}$).

\begin{figure}[h!]
\begin{subfigure}{.5\textwidth}
    \centering
    \includegraphics[width=\textwidth]{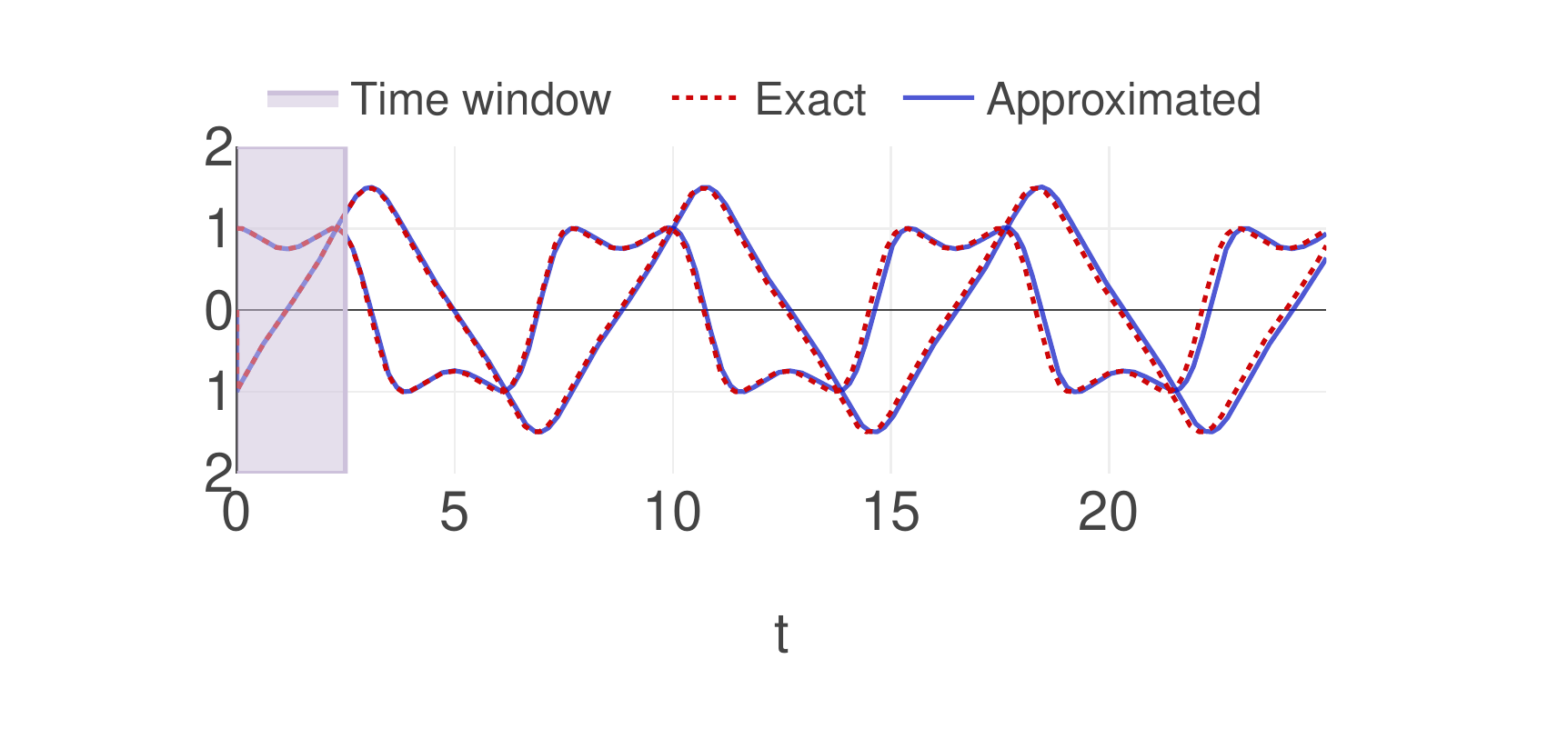}
    \caption{$a=1.12$}
    \label{fig:ODE1}
\end{subfigure}
\begin{subfigure}{.5\textwidth}
    \centering
    \includegraphics[width=\textwidth]{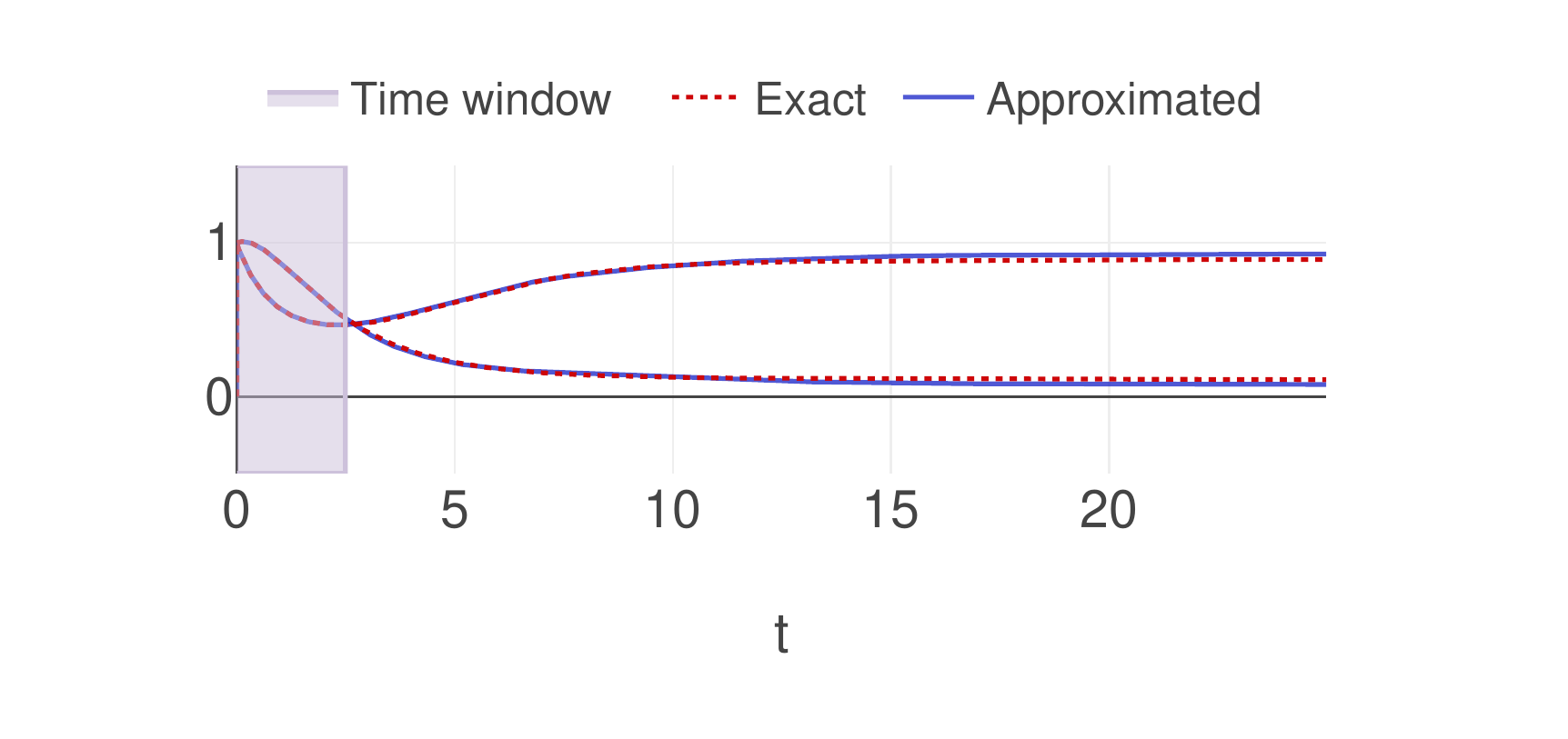}
    \caption{$a=1.12$}
    \label{fig:ODE2}
\end{subfigure}
\begin{subfigure}{.5\textwidth}
    \centering
    \includegraphics[width=\textwidth]{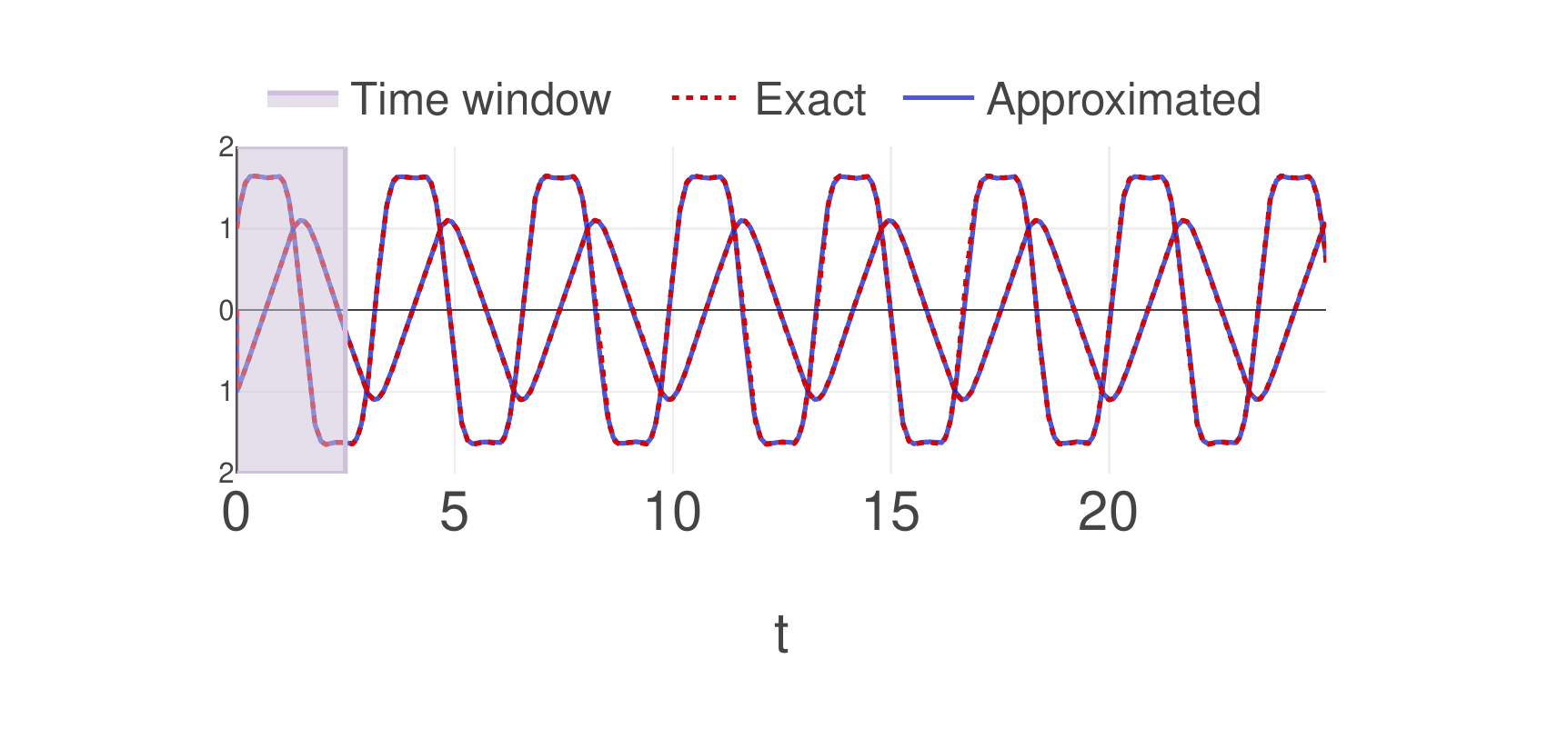}
    \caption{$a=5.22$}
    \label{fig:ODE3}
\end{subfigure}
\begin{subfigure}{.5\textwidth}
    \centering
    \includegraphics[width=\textwidth]{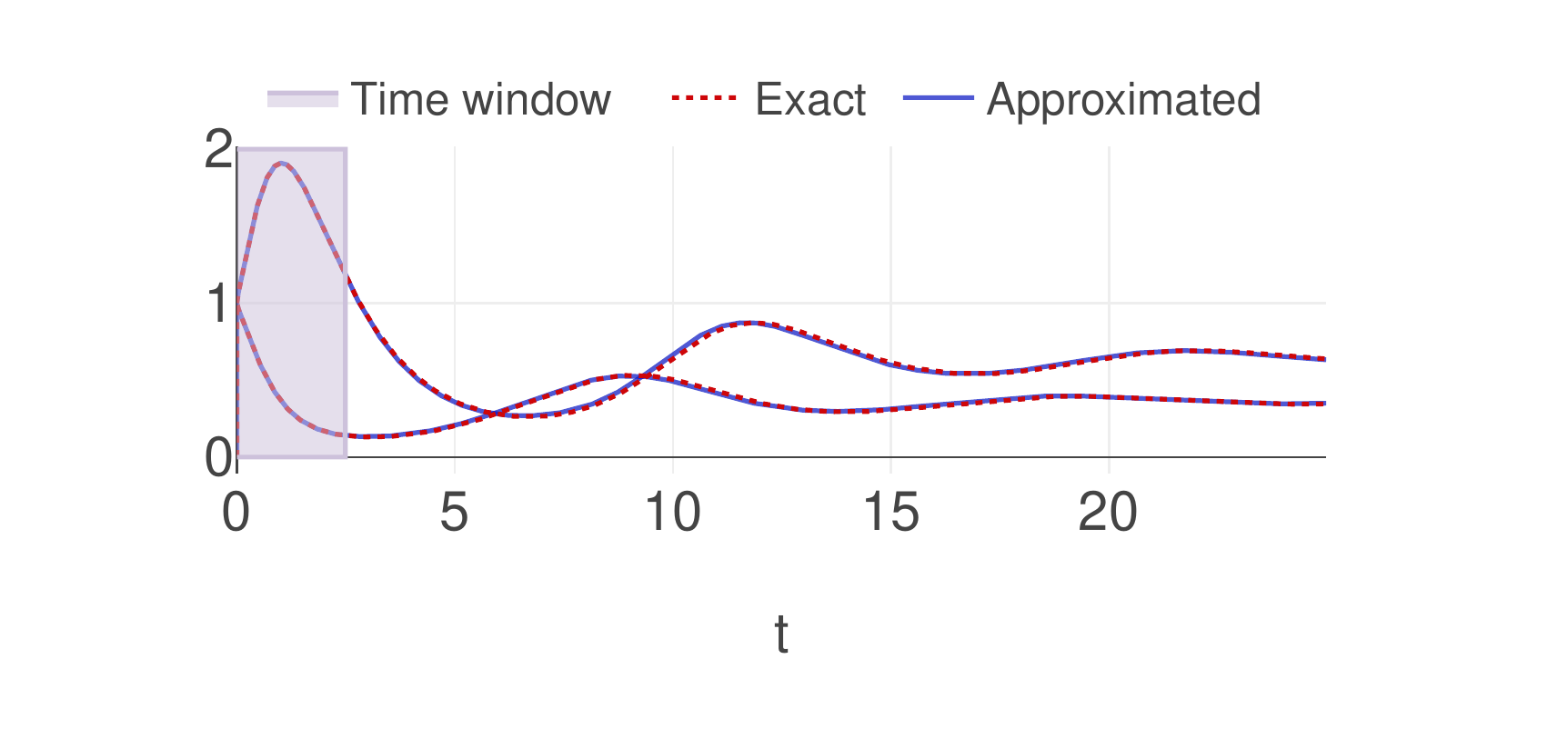}
    \caption{$a=2.78$}
    \label{fig:ODE4}
\end{subfigure}
\begin{subfigure}{.5\textwidth}
    \centering
    \includegraphics[width=\textwidth]{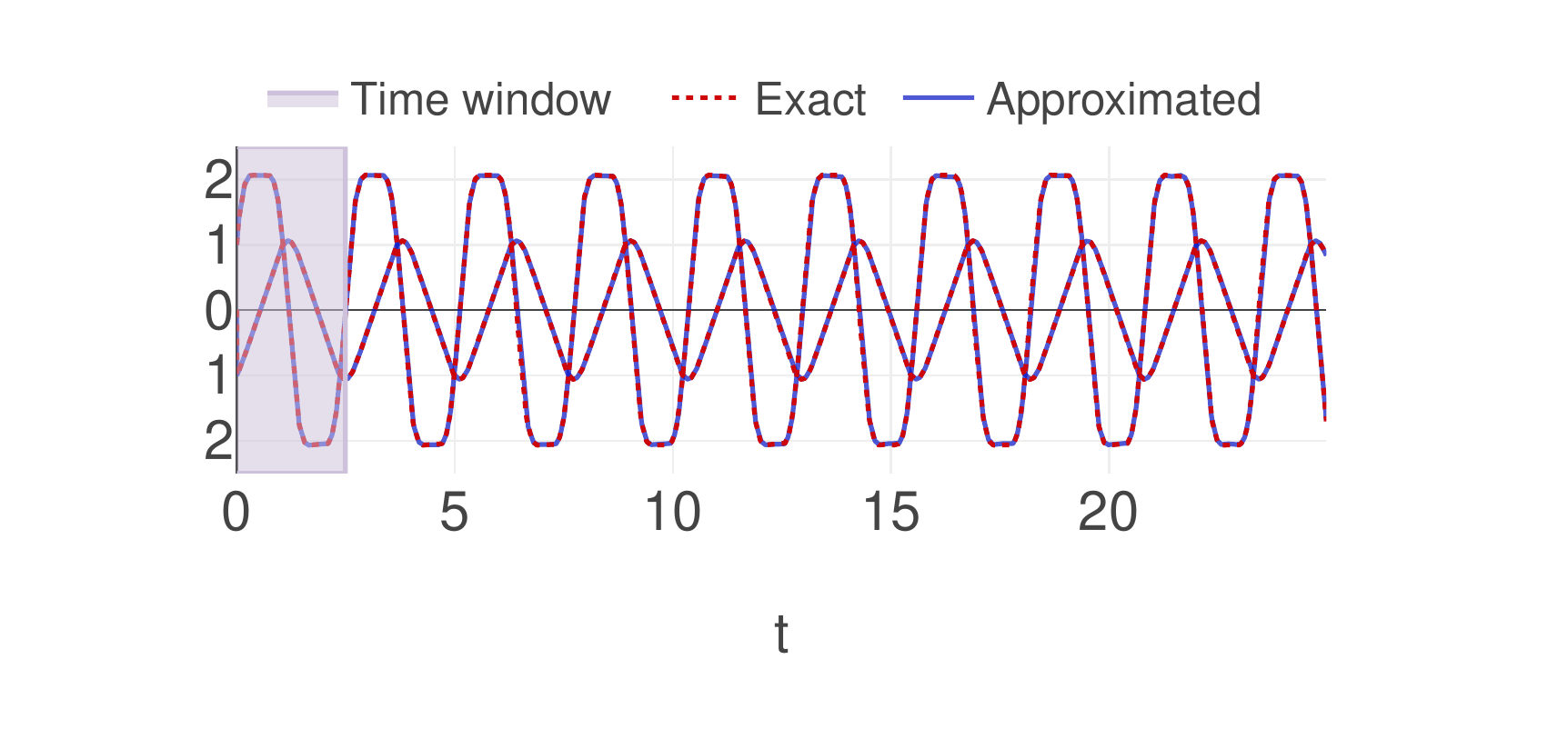}
    \caption{$a=8.43$}
    \label{fig:ODE5}
\end{subfigure}
\begin{subfigure}{.5\textwidth}
    \centering
    \includegraphics[width=\textwidth]{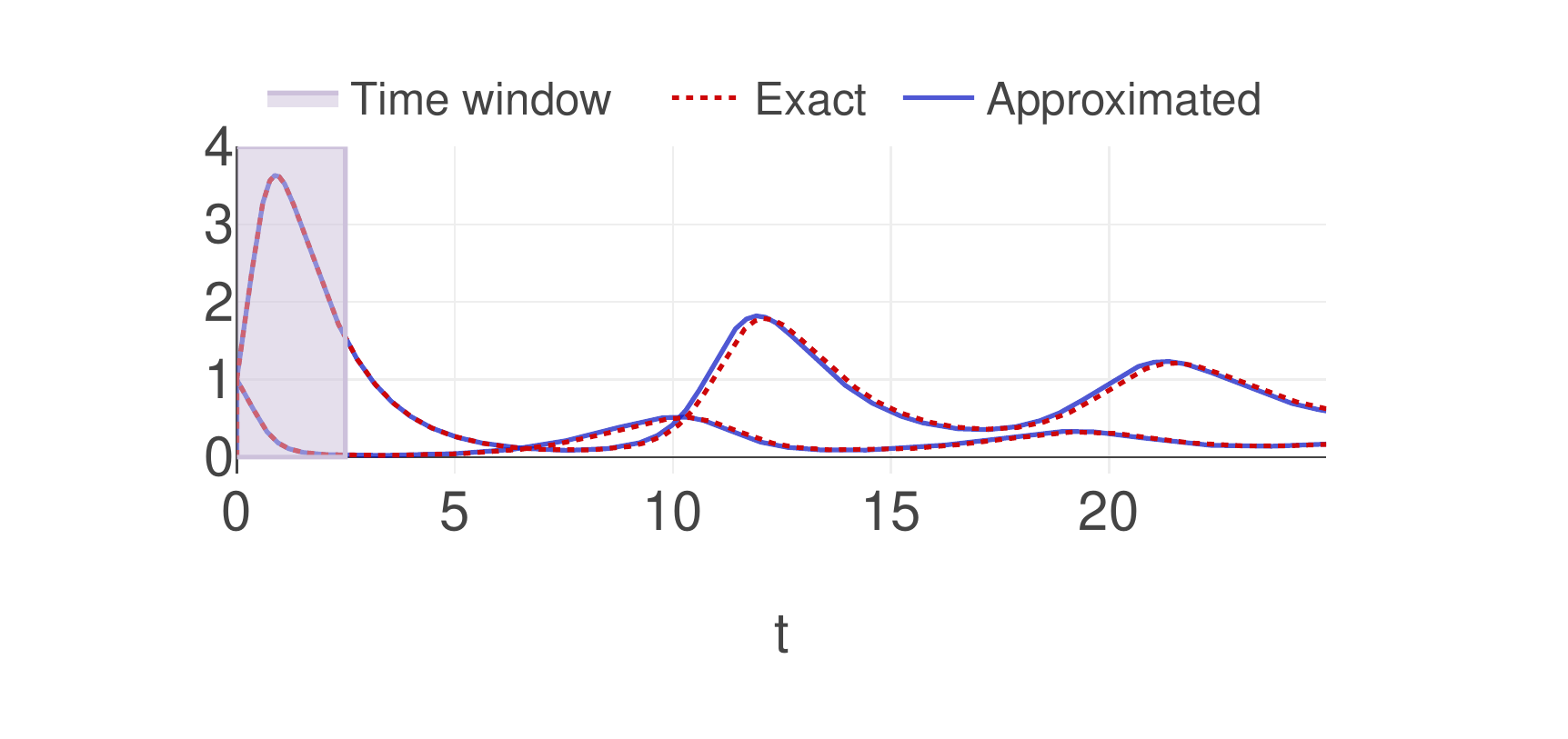}
    \caption{$a=4.87$}
    \label{fig:ODE6}
\end{subfigure}
\caption{ \small On the left column, the exact and the predicted evolution of the Duffing Oscillator system for different values of the parameter $a$ (\ref{fig:ODE1},\ref{fig:ODE2},\ref{fig:ODE5}). On the right one, the exact and the predicted evolution of the non-linear Predator Prey system (\ref{fig:ODE2},\ref{fig:ODE4},\ref{fig:ODE6}).}
\label{fig:ODE}
\end{figure}

As we can in general see from Figure \ref{fig:ODE}, the architecture succeeds in reproducing the systems' dynamics also for parameters not included in the training set, thus it is able to generalize the system's parametric behaviour. 

Furthermore, in order to investigate the role of some architectural parameters, tests have been conducted on the influence that the number of time-steps predicted per iteration, $m$, have in the accuracy of the predictions. 

\begin{figure}[h!]
    \centering
    \includegraphics[width=.7\textwidth]{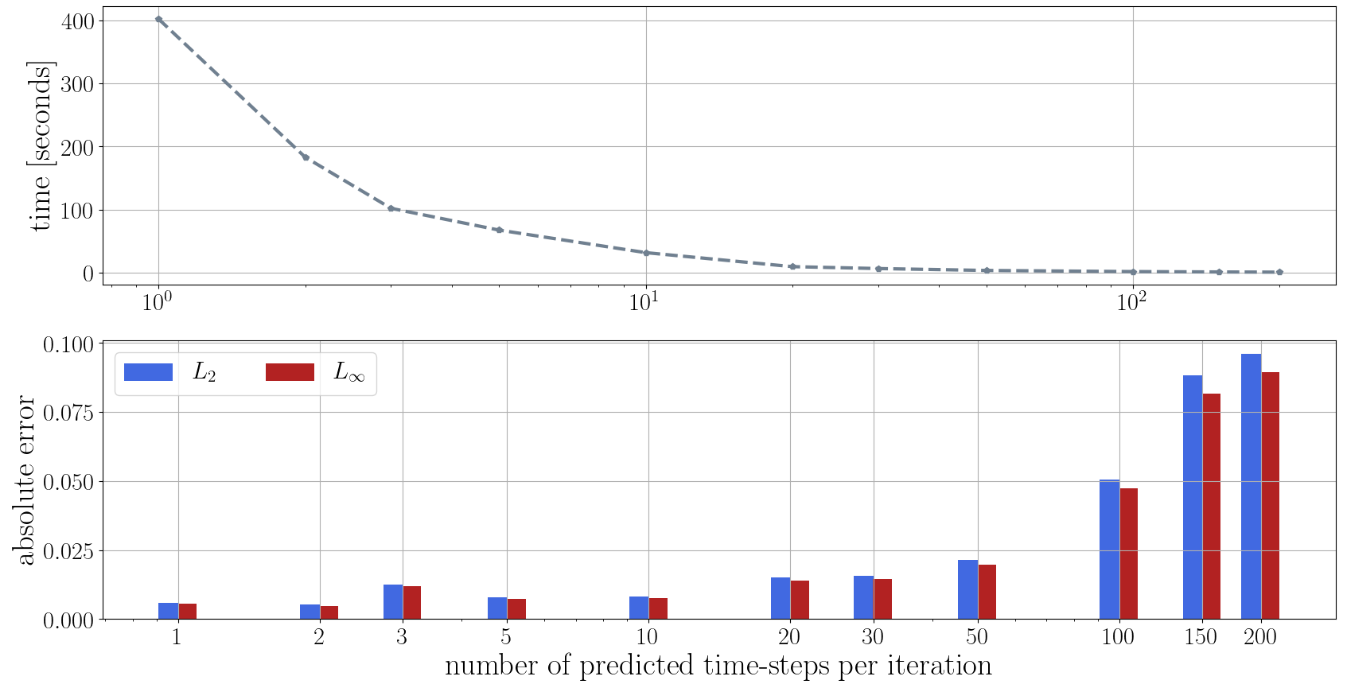}
    \caption{Time needed for the predictions and mean absolute error committed in the first 1000 time-steps versus the value imposed on $m$. These tests have been performed on both the ODEs systems previously described for $24$ testing parameters sampled in the respective parameter spaces.}
    \label{fig:t-err}
\end{figure}

As it can be seen in Figure \ref{fig:t-err}, enlarging $m$ could bring some advantages in terms of the amount of time needed to predict a certain number $T$ of time-steps (less iteration cycles are required). On the other hand, the drawback for large values of $m$ appears to be the reducing of generalization capabilities of the architecture.

Nevertheless, the size of such models is too small to observe any speed-up, therefore these examples only want to serve as introductory analysis. Indeed, referring to the python package \emph{tfdiffeq}\footnote{https://github.com/titu1994/tfdiffeq/tree/master/tfdiffeq}, we take as our baseline the time spent by its function \emph{odeint()} for the integration of the system. Those times for the above ODEs amount to $3 s$ on average (with the default available solver, that implements an adaptive Runge-Kutta algorithm). On the other hand, observing the graphics in Figure \ref{fig:t-err}, we find that the two better alternatives with our framework are obtained with $m=5$ or $m=10$ (for what concerns a trade-off between accuracy and time needed). The correspondent prediction-times are respectively $70 s$ and $31 s$, both of which under-perform our baseline by one order of magnitude.

\begin{figure}[h!]
    \centering
    \includegraphics[width=.7\textwidth]{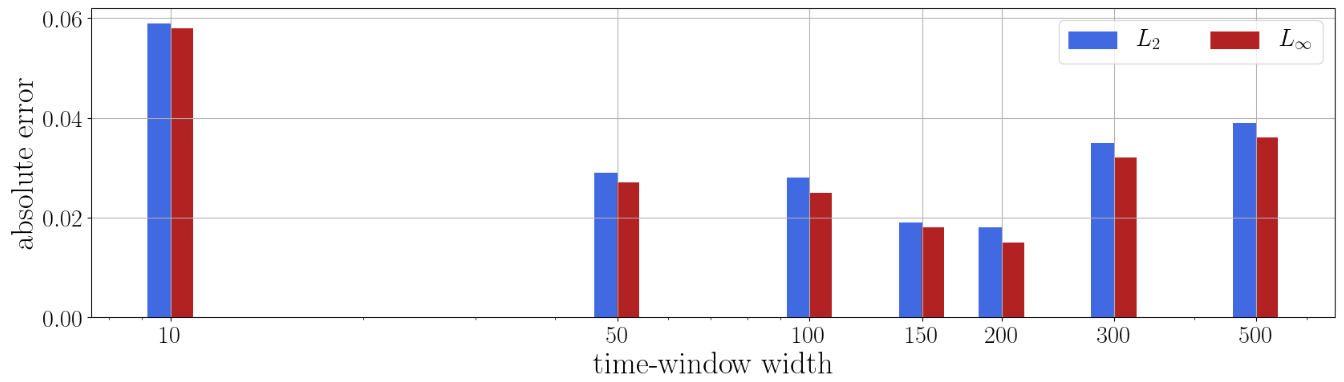}
    \caption{mean absolute error committed in the first 1000 time-steps versus the value imposed on $w$. These tests have been performed on both the ODEs systems previously described for $24$ testing parameters sampled in the respective parameter spaces.}
    \label{fig:w-err}
\end{figure}

Moreover, tests have been also performed on the correlation between the time-window size and the accuracy of the predictions. As is can be seen observing Figure \ref{fig:w-err}, an initial decreasing trend in the relative error appears evident with the increase of the window size. This could be justified by the fact that a larger time-window implies more exact information at the beginning of the prediction iteration, thus bringing a slower error propagation in the process. We can also note that increasing too much $w$ actually does not bring any additional improvement on the error (in this case from $w=200$ on). The time measurements are here not reported, as no difference is  encountered varying the time-window width. 

\section{Rayleigh-B\'{e}nard cavity flow}\label{section:cavity}

In order to extend our tests to larger systems, we present here the Rayleigh-B\'{e}nard cavity flow: a benchmark example that has been introduced in \cite{Roux:GAMM} and widely used since then, for example in \cite{Gelfgat:Ref11}, \cite{PR15} and \cite{HessQuainiRozza_ACOM_2022}. It considers the incompressible Navier-Stokes equations in a rectangular cavity. In particular, the model describes an important process in  semiconductor crystal growth \cite{KAKIMOTO1995191}, as it models the flow in the molten semiconductor material.

\subsection{Model description}

The incompressible Navier-Stokes equations describe viscous, Newtonian flow 
in the computational domain $\Omega \subset \mathbb{R}^d$. 
The unknowns are the vector field velocity $\bm u$ and scalar field pressure $p$. The incompressible Navier-Stokes equations are given as
\begin{eqnarray}\label{NS-1}
\frac{\partial {\bm u}}{\partial t}+({\bm u}\cdot  \nabla {\bm u})-\nu\Delta {\bm u}+\nabla p={\bm f} &\qquad\mbox{in } \Omega\times(0,T],\\ \label{NS-2}
\nabla \cdot {\bm u}=0&\qquad\mbox{in } \Omega\times (0,T].
\end{eqnarray}

\noindent where the kinematic viscosity is denoted $\nu$, the body forcing $\bm f$, and time interval as $T$. 
The spatial dimension is either $d=2$ or $d=3$, while boundary and initial conditions are provided as
\begin{eqnarray}
{\bm u}={\bm u}_0&\qquad&\mbox{in } \Omega \times \{0\}, \label{IC}\\
{\bm u}={\bm u}_D&\qquad&\mbox{on } \partial\Omega_D\times (0,T], \label{BC-D} \\
-p \bm n + \nu \frac{\partial {\bm u}}{\partial \bm n}= \bm g &\qquad&\mbox{on } \partial\Omega_N\times, (0,T], \label{BC-N}
\end{eqnarray}

where $\partial\Omega_D \cap \partial\Omega_N = \emptyset$ and $\overline{\partial\Omega_D} \cup \overline{\partial\Omega_N} =  \overline{\partial\Omega}$. Here, ${\bm u}_0$, ${\bm u}_D$, and $\bm g$ are given and $\bm n$ denotes the outward pointing unit normal vector on the boundary $\partial\Omega_N$. The boundary $\partial \Omega_D$ is called the Dirichlet boundary and $\partial \Omega_N$ the Neumann boundary.

Let $L^2(\Omega)$ denote the space of square integrable functions in $\Omega$ and $H^1(\Omega)$ the space of functions belonging to $L^2(\Omega)$ with weak first derivatives in $L^2(\Omega)$. 
Define the sets 
\begin{eqnarray}
{\bm V} &:=& \left\{ {\bm v} \in [H^1(\Omega)]^d:  {\bm v} = {\bm u}_D \mbox{ on }\partial\Omega_D \right\},  \\
{\bm V}_0 &:=& \left\{{\bm v} \in [H^1(\Omega)]^d:  {\bm v} = \boldsymbol{0} \mbox{ on }\partial\Omega_D \right\}. 
\end{eqnarray}
The variational form of \eqref{NS-1}--\eqref{BC-N} is given by: find $({\bm u},p)\in {\bm V} \times L^2(\Omega)$, with ${\bm u}$ satisfying the initial condition \eqref{IC}, such that
\begin{align}
&\int_{\Omega} \frac{\partial{\bm u}}{\partial t}\cdot{\bm v} \mathbf{x}+\int_{\Omega}\left({\bm u}\cdot\nabla {\bm u}\right)\cdot{\bm v} \mathbf{x}
+\nu \int_{\Omega}\nabla {\bm u}\cdot \nabla{\bm v} \mathbf{x}- \int_{\Omega}p\nabla \cdot{\bm v} \mathbf{x} \nonumber \\
&\hspace{3cm} =\int_{\Omega}{\bf f}\cdot{\bm v} \mathbf{x}  + \int_{\partial \Omega_N}{\bf g}\cdot{\bm v} \mathbf{x} ,
\qquad\forall\,{\bm v} \in {\bm V}_0, \label{eq:weakNS-1}\\
& \int_{\Omega}q\nabla \cdot{\bm u} \mathbf{x} =0, \qquad\forall\, q \in L^2(\Omega).  \label{eq:weakNS-2}
\end{align}
Consider as computational domain $\Omega$ the rectangle with aspect ratio 4, i.e., a rectangle of height 1 and length 4.
The whole boundary is a no-slip boundary, so that $\partial\Omega_D = \partial\Omega$ and $\bm u_D = 0$.
The body forcing ${\bm f}$ is given by
\begin{align}
    {\bm f} = (0, \text{Gr}\nu^2 x)^T, \label{eq_f_cavity}
\end{align}
\noindent where $x$ is the horizontal coordinate and Gr is the Grashof number.
The Grashof number is a dimensionless number that describes the ratio of the buoyancy to viscous forces. 
\subsection{Discretization}

\begin{figure}[h!]
    \centering
    \includegraphics[scale=.4]{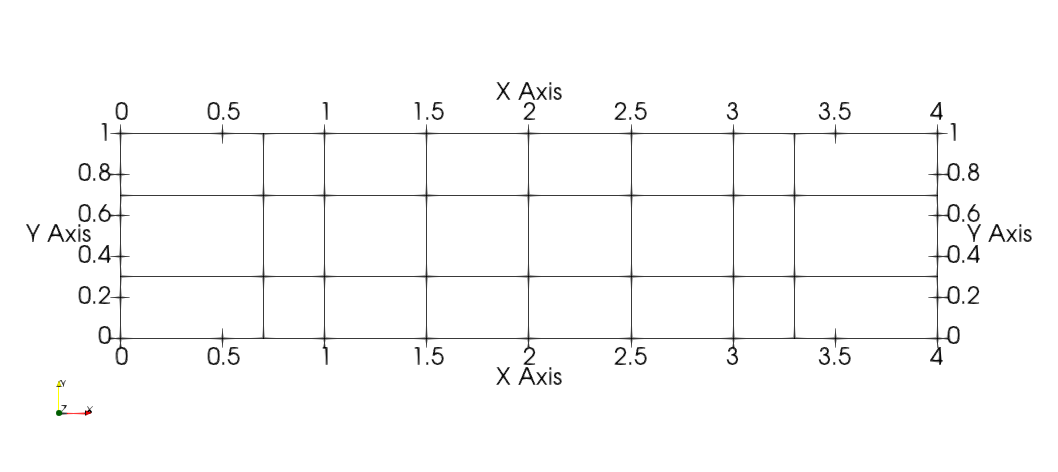}
    \caption{\small The computational mesh of the cavity is composed of 24 rectangles.}
    \label{fig:mesh_cavity}
\end{figure}

The numerical discretization method employed is the spectral element method
 \cite{karniadakis1999spectral}, which uses high-order polynomial ansatz functions over a coarse mesh, see Fig.~\ref{fig:mesh_cavity}.
The time-stepping scheme is an IMEX scheme of order 2 (IMplicit-EXplicit, see
\cite{Guermond_Shen_VCS}, \cite{Karniadakis_Orszag_Israeli_splitting_methods}), which is a standard option of the used PDE solver \emph{Nektar++}\footnote{https://www.nektar.info/}.

Our numerical studies will focus on the parameter domain $\Theta = [100 \cdot 10^3, 150 \cdot 10^3]$, and a full-order solution is computed at $Gr = 150\cdot 10^3$ over a long time interval to ensure that the limit cycle is reached. Then, each solution of interest in the interval $[100 \cdot 10^3, 150 \cdot 10^3]$ is initialized with the solution at $Gr = 150\cdot 10^3$. 

The time step is set to $1 \cdot 10^{-6}$ and for the tests $5 \cdot 10^{5}$ time steps have been computed.


\subsection{Model order reduction}\label{POD}

Our Model Order Reduction technique aims to reduce
the cost of the full order solution computation by breaking it into two parts: a computationally expensive \emph{offline phase}, and a computationally efficient \emph{online phase}.

Indeed, the offline phase is the most time consuming because it comprehends both the collection of the full order solutions, and the training of the two stages architecture with these. On the other hand, the online phase is intended to be particularly fast, as it consists only in the computation of the first exact window $W_{new}$ and in the iterative time-step prediction by the framework (see Section \ref{architecture}).

Recalling what has been already analyzed in the ODE case in section~\ref{section:ODE}, the number of time-steps $w$ belonging to the time-window has an influence on the accuracy of the predictions, as a lower one, in general, implies a stronger error propagation. Despite of this, in the case of high-dimensional models, the computation of a long initial exact window can be as time consuming as the iterative prediction phase. Therefore, a compromise has to be made.

Moreover, to further reduce the order of the model, previous operations are performed on the full order solutions. 

In the first place, the velocity field solutions at every time step in the time interval of interest $T$, which are real vectors of high dimension $N$ ($N$
referring to size of the spatial discretization), are projected on a lower dimensional space through the Proper Orthogonal Decomposition approach (POD). It consists in finding a certain number of POD modes that reduces the dimension of the \emph{snaphots matrix} $\mathbf{S} = [\mathbf{s}_1,\dots,\mathbf{s}_T]$ (being $[\mathbf{s}_1,\dots,\mathbf{s}_T]$ the $T$ $N$-dimensional full order solutions). Such basis are computed through the Singular Value Decomposition of $\mathbf{S}$: $$\mathbf{S} = \mathbf{U}\mathbf{\Sigma}\mathbf{V^T},$$ where the columns of the unitary matrix $\mathbf{U}$ are the POD modes, and the diagonal matrix $\mathbf{\Sigma}$ contains the corresponding singular values in decreasing order. Considering only the first $N_{POD}$ rows of $\mathbf{U}$ according to the error we are willing to commit, we obtain the reduced order representation: $$\mathbf{S_{POD}} =  \mathbf{U_{N_{POD}}^T}\mathbf{S}.$$

In second place, because of the very small time-step required by the cavity simulations, in the collection of the training data-set for the two stages architecture only $1$ every $100$ time-steps is considered, so that also the online phase is accelerated having a larger prediction step.


\section{Results} \label{sec:res}

As discussed in Section \ref{POD}, the first step of our model order reduction approach in the case of significantly large systems is Proper Orthogonal Decomposition. The \emph{snapshots matrix} $\mathbf{S}$ is thus composed by full-order solutions, which have been collected for $6$ equispaced values of Grashof number: $$\{\mathbf{\theta}_{t}^i\}_{i=1}^n = \{100\cdot 10^3,110\cdot 10^3,120\cdot 10^3,130\cdot 10^3,140\cdot 10^3,150\cdot 10^3\}.$$ SVD is then performed to obtain $N_{POD}$ basis, chosen to achieve a certain level of accuracy in the approximation.

Indeed, the accuracy can be derived as the ratio between the summation of the singular values correspondent to the considered POD basis and the summation of the whole $\Sigma$ diagonal. In our case, a number of $140$ and $147$ POD basis is identified to achieve a $99.99\%$ level of accuracy, respectively for the horizontal and vertical dimension. 

As a preliminary analysis, we are interested in the number of POD modes coefficient it is actually convenient to consider during the training in the offline phase. Indeed, even if the more modes are involved, the more POD succeeds in an accurate decomposition, the architecture could encounter more difficulties in predicting a larger number of outputs rather then a lower one.

\begin{figure}[h!]
    \centering
    \includegraphics[width=.85\textwidth]{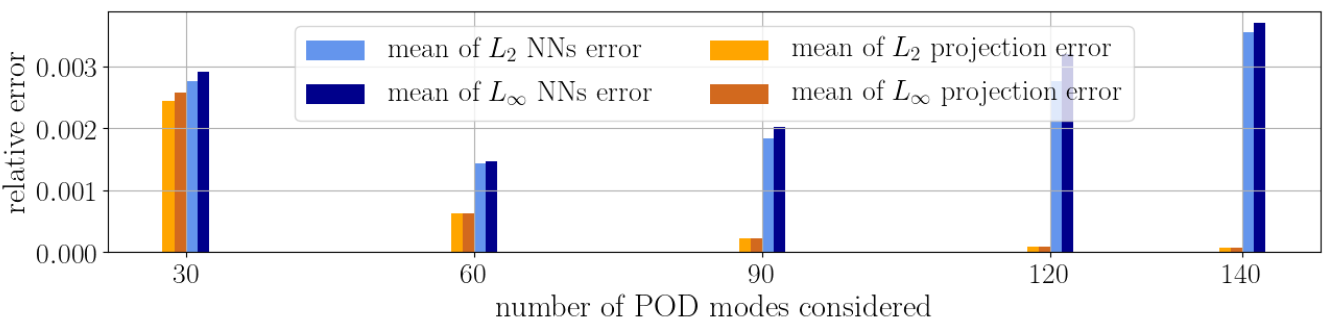}
    \caption{Projection mean error and NNs mean error with respect to different numbers of POD coefficients taken into consideration.}
    \label{POD_modes_study}
\end{figure}

For this reason, different training phases have been firstly performed considering each time a different number of POD modes, and the results are shown in Figure \ref{POD_modes_study}. A distinction is made between the \emph{projection error}, \emph{i.e.} the error due to POD, and the \emph{NNs error}, which is the prediction error of the architecture. This latter one clearly depends both on the projection error and on the generalization capabilities of the framework. 

In particular, the chosen framework's parameters are $m=6$ and $w=100$, correspondent to a time-window of $10000$ real time-steps and to an horizon of $600$ future time-steps predicted per iteration.

The results displayed in Figure \ref{POD_modes_study} are obtained with $10$ testing parameters sampled in the Grashof number space. Starting from an initial exact window, the evolution in time of the correspondent POD coefficients is computed for the next $500000$ real time-steps. 

From Figure \ref{POD_modes_study} it is clear that, while the projection error monotonously decreases, the NNs error reaches a minimum and then increases with the number of POD coefficients considered. Owing to such evidence, the following results are computed considering only the first $60$ POD coefficients.

\begin{figure}[h!]
    \centering
    \includegraphics[width=.85\textwidth]{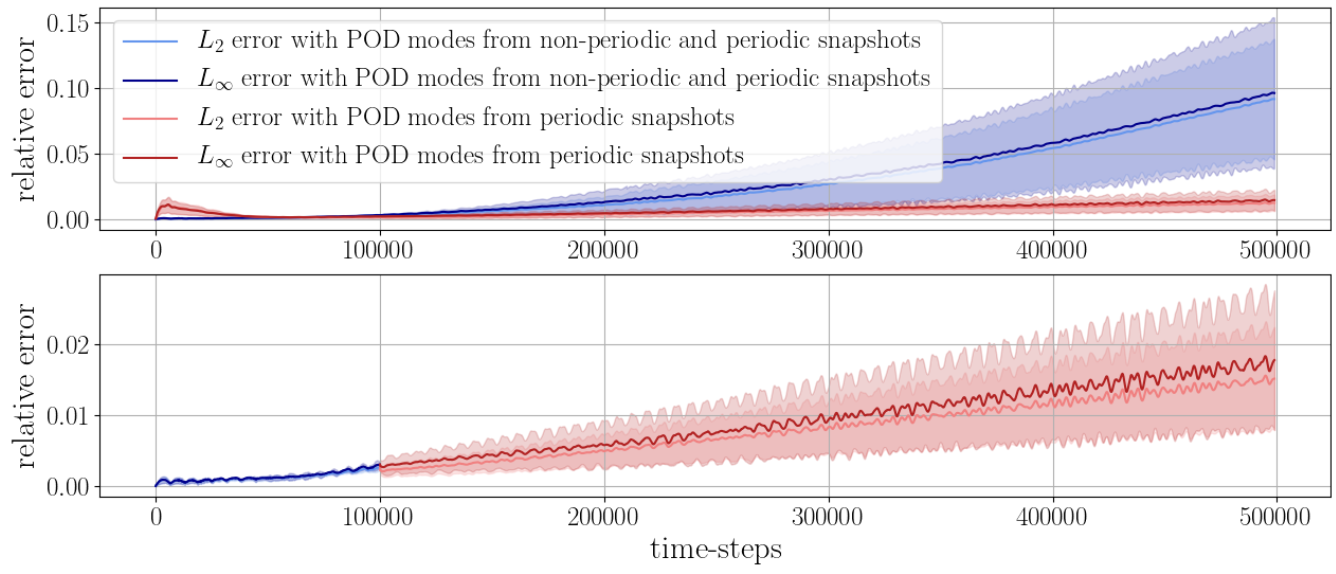}
    \caption{Mean relative error propagation in time. In the upper graphic the two cases of POD: the blue one considering both the non-periodic and periodic snapshots and the red one taking only the periodic ones. In the lower graphic the error propagation in the case of pipelining in time the two PODs is displayed.}
    \label{fig:cavity}
\end{figure}

Considering thus the $10$ testing values of the Grashof number randomly sampled, we report the relative error evolution in Figure \ref{fig:cavity} with the blue lines in the upper graphic. As it can be seen, for the represented $5000$ time-steps, the accuracy tends to decrease quickly, touching mean error peaks of $10\%$.

\begin{figure}[h!]
\begin{subfigure}{\textwidth}
    \centering
    \includegraphics[width=\textwidth]{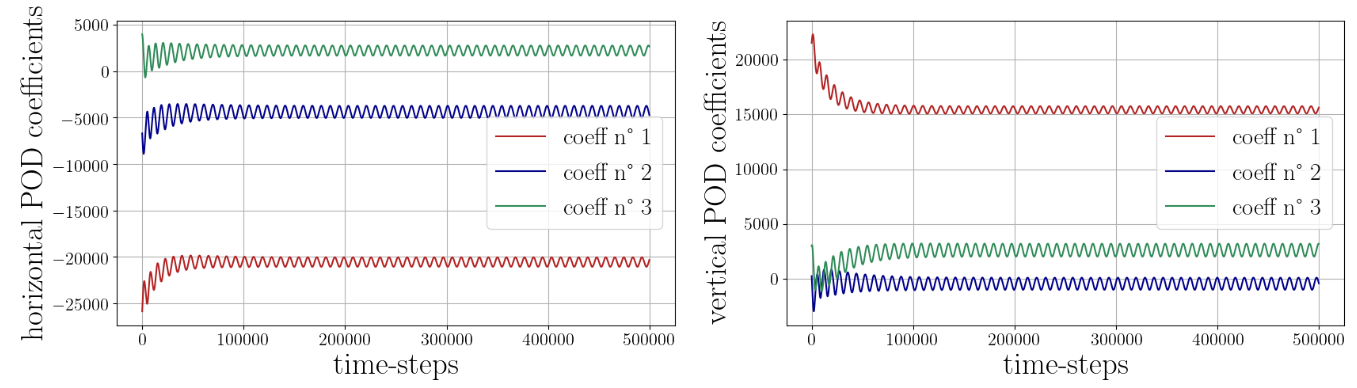}
\end{subfigure}
\caption{Evolution in time of the first 3 POD basis' coefficients correspondent to the exact solution of the Cavity problem with Grashof number equal to $100000$. Being the problem treated in two dimensions, we find on the left the coefficients related to the horizontal axis, and on the right the ones for the vertical axis.}
\label{fig:modes}
\end{figure}

Investigating more the problem, we note, from Figure \ref{fig:modes}, that two different behaviours in the snapshots time evolution can be observed: an initial \emph{swing-in} non-periodic phase, and a successive periodic one. Given that, we try to further reduce the dimensionality by searching new POD basis considering only the periodic-part's snapshots.

In this way, we find that the number of POD basis needed to reach the same purpose as before ($99.99\%$ of accuracy) is significantly lower: respectively $37$ and $41$ for the horizontal and vertical axis. Therefore, wanting  to investigate the error committed in the prediction of such $37$ modes coefficients, we train the two stages architecture on the new collected training set. 

Performing the same tests with $10$ different values of the Grashof number, we obtain that the framework is now able to approximate the velocity field with a much lower relative error, as it can be seen in Figure \ref{fig:cavity} in the upper graphic with the red lines. On the other hand, with this second approach the initial non-periodic part approximation, well predicted in the first case, results generally worsened.

Hence, not considering the swing-in phase time-steps in the \emph{snapshots matrix}, the number of POD basis $N_{POD}$ needed to achieve a certain accuracy of the approximation considerably decreases, as a lower number of singular values turn out to be significantly energetic. As a consequence, error propagation seems to be attenuated in the second case, owing to the lower number of different coefficients to be predicted. On the other hand, with less significant modes the initial swing-in phase seems to not be sufficiently well approximated.

The solution we propose, displayed in the lower graphic of Figure \ref{fig:cavity}, is to form a pipeline with the two architectures previously trained. More precisely, we let the first framework predict the first $N_I$ time-steps correspondent to the non-periodic behaviour, and then a basis change is performed to project the last $w$ approximated time-steps $\mathbf{S_{POD_1}}$ into the lower dimensional space (from $140$ to $37$ dimensions). At this point, the new input window $\mathbf{S_{POD_2}}$ is obtained for the second framework, which can now deal with the approximation of the periodic part. In particular, the matrix for basis change is obtained as: $$\mathbf{M} = \mathbf{U_{N_{POD_2}}^T}\cdot \mathbf{U_{N_{POD_1}}},$$ and the new input for the second framework is: $$\mathbf{S_{POD_2}} = \mathbf{M}\cdot \mathbf{S_{POD_1}}.$$

\begin{figure}[h!]
    \centering
    \includegraphics[width=.85\textwidth]{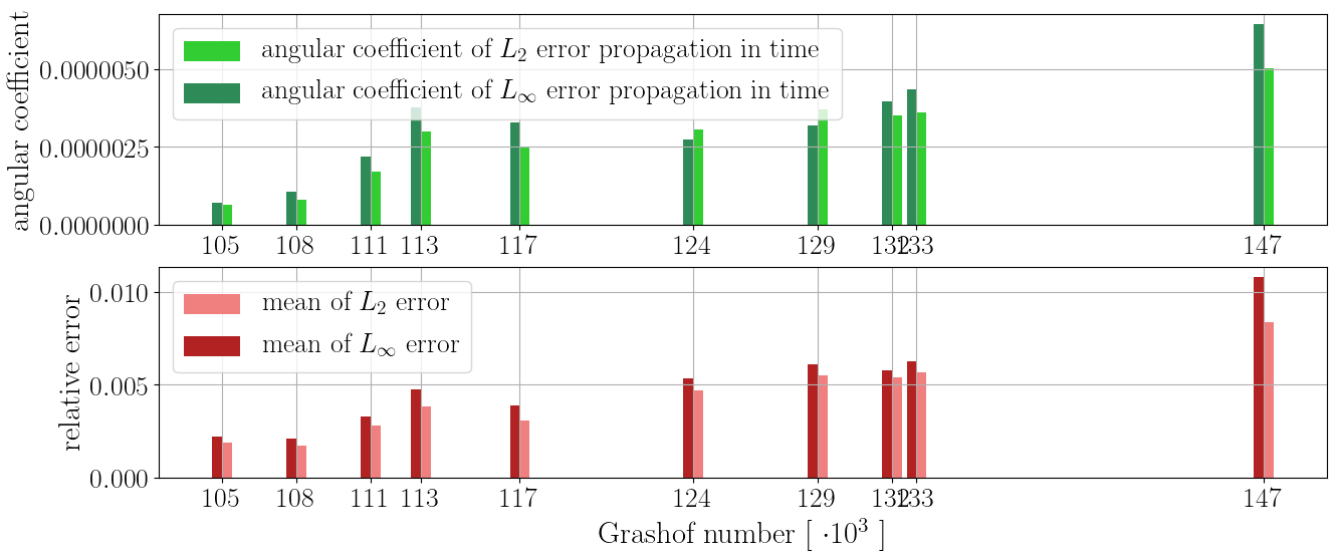}
    \caption{Mean error and correspondent angular coefficient of its linear regression in time computed in the testing cases taken into consideration.}
    \label{fig:evolution_error}
\end{figure}

As we can observe in the comparison between the upper and lower graphic in Figure \ref{fig:cavity}, the relative error's evolution for what concerns the periodic part does not significantly differ in the two cases, thus this pipeline is effective and does not generate a worse error propagation. The testing results are then reported in Figure \ref{fig:evolution_error} in terms of mean errors on a time horizon of $500000$ time-steps, and of the angular coefficients correspondent to their linear regression.

In general, the big advantage that the application of the two stages architecture brings, involves the time needed to obtain a new Grashof number's velocity field solution. Indeed, the computational time for predicting $500000$ real time-steps is reduced to $5$ minutes on average, versus the $3$ hours spent with the \emph{Nektar solver} to obtain the high-dimension evolution of a new Grashof number's solution.

\begin{figure}[h!]
\centering
\includegraphics[width=1\textwidth]{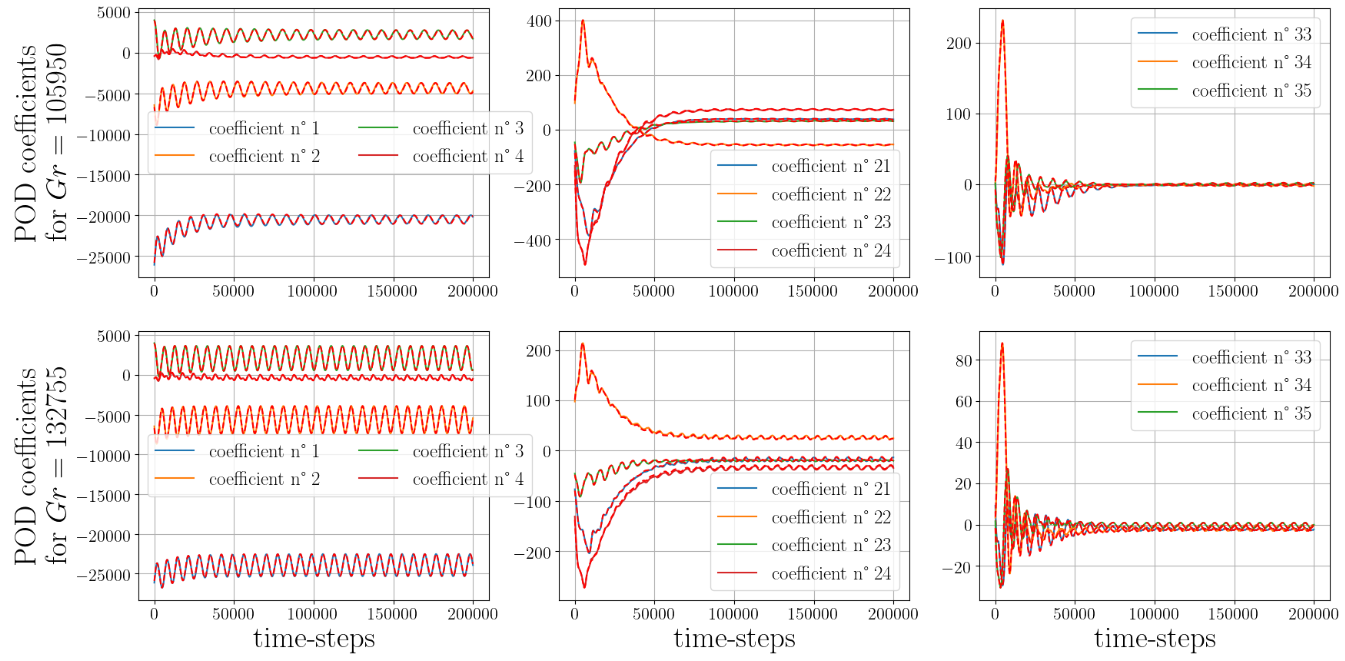}
\caption{Example of some modes coefficients correspondent to $Gr=132775$, for the horizontal and vertical axis respectively in the first and second row. In coloured solid lines the predicted time-evolution, while in red dashed-lines the exact one. The coefficients here reported are related to the $37$ and $41$ POD basis discussed above.}
\label{fig:modes_comp}
\end{figure}

\begin{figure}[h!]
\centering
\includegraphics[width=1\textwidth]{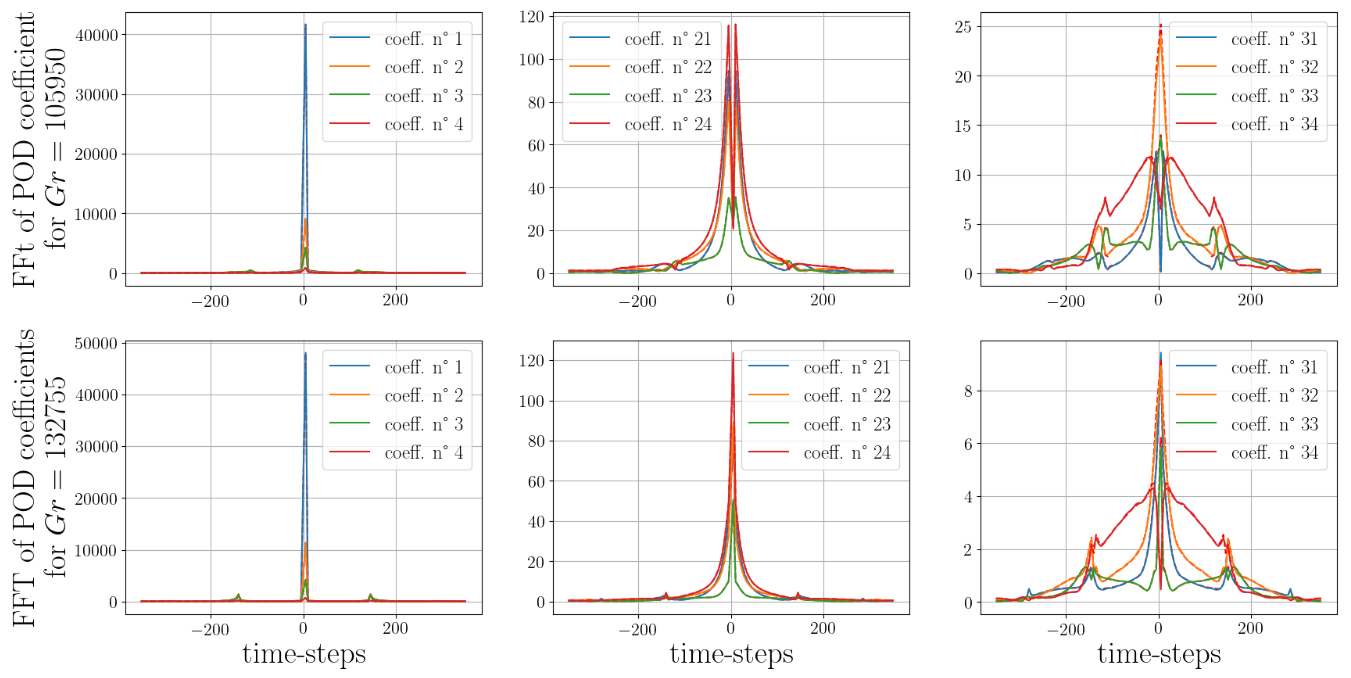}
\caption{\footnotesize Example of the Fourier Transform of some modes coefficients for the horizontal and vertical axis respectively in the first and second row ($Gr=132775$). In coloured solid lines the predicted time-evolution Fourier Transform, while in red dashed-lines the exact one.}
\label{fig:FFT_modes_comp}
\end{figure}

\begin{figure}[h!]
\begin{subfigure}{.5\textwidth}
    \centering
    \includegraphics[width=\textwidth]{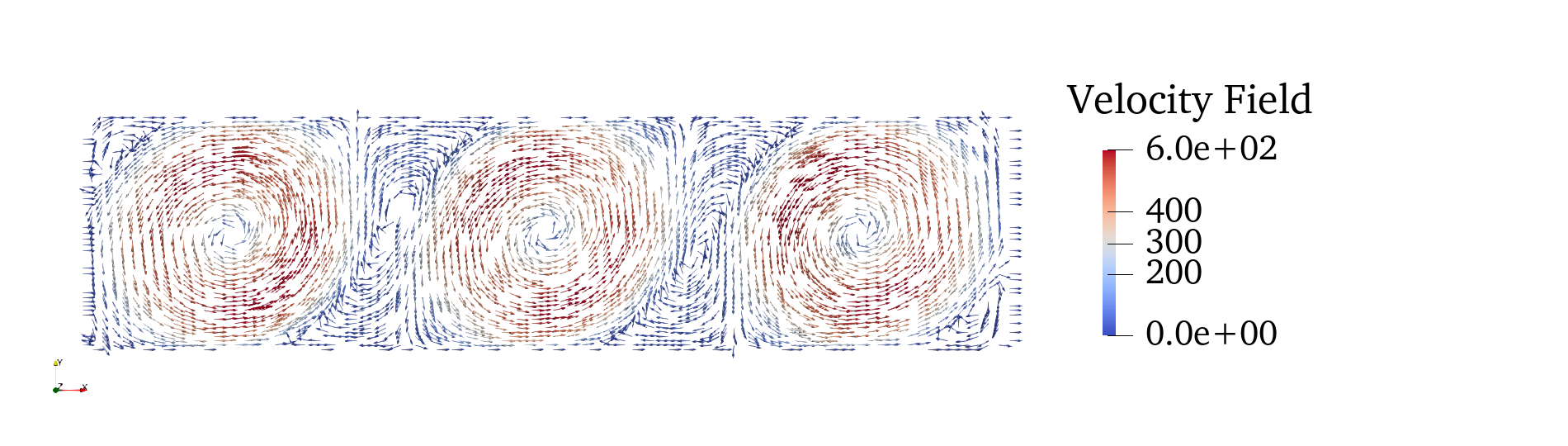}
\end{subfigure}
\begin{subfigure}{.5\textwidth}
    \centering
    \includegraphics[width=\textwidth]{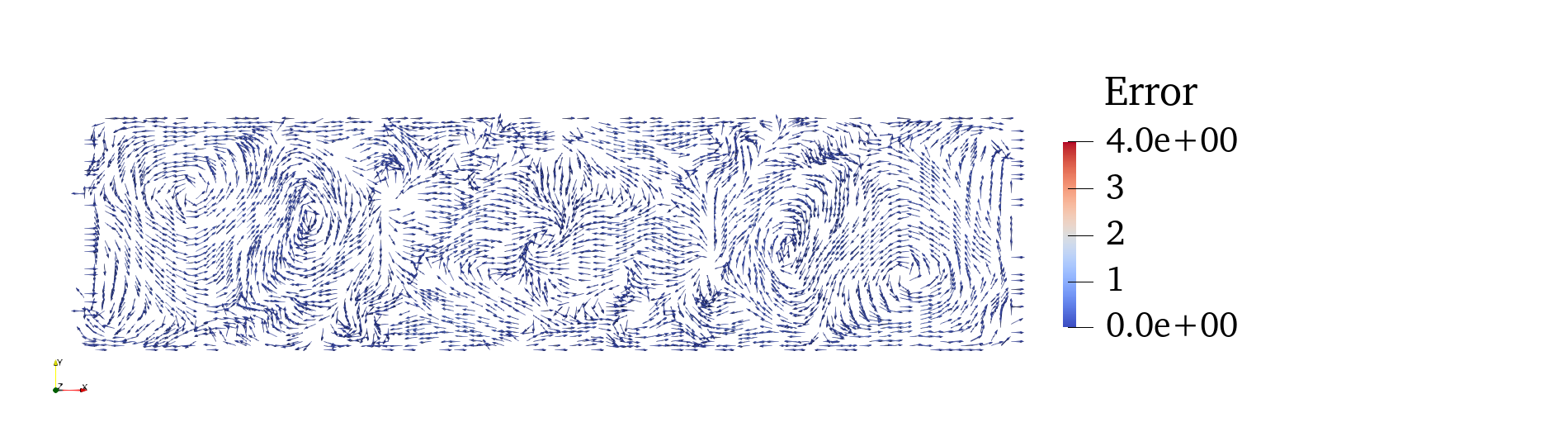}
\end{subfigure}
\begin{subfigure}{.5\textwidth}
    \centering
    \includegraphics[width=\textwidth]{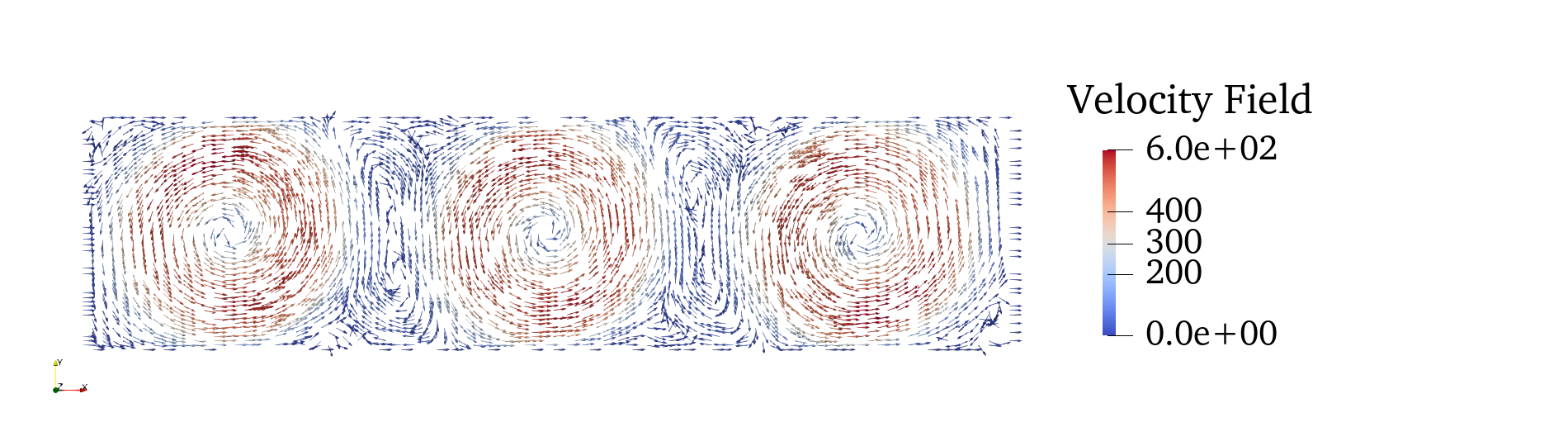}
\end{subfigure}
\begin{subfigure}{.5\textwidth}
    \centering
    \includegraphics[width=\textwidth]{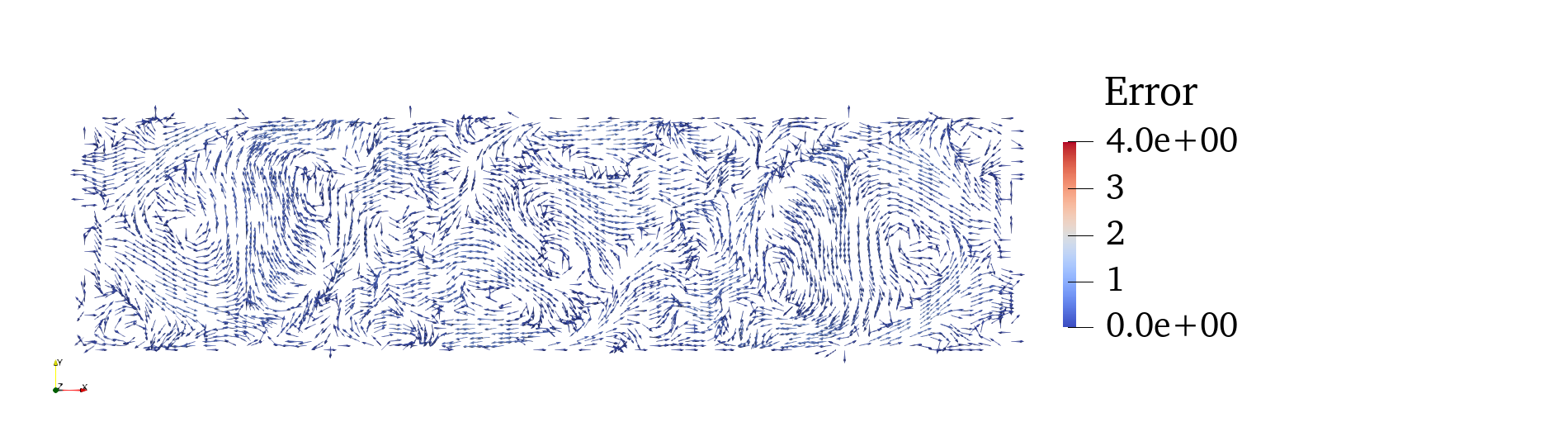}
\end{subfigure}
\begin{subfigure}{.5\textwidth}
    \centering
    \includegraphics[width=\textwidth]{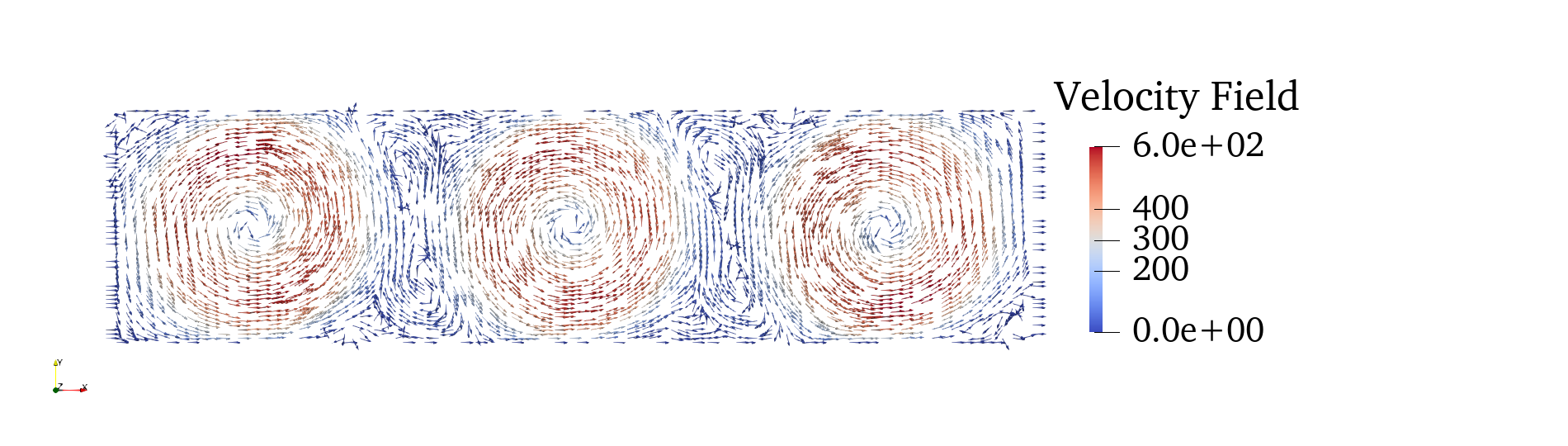}
\end{subfigure}
\begin{subfigure}{.5\textwidth}
    \centering
    \includegraphics[width=\textwidth]{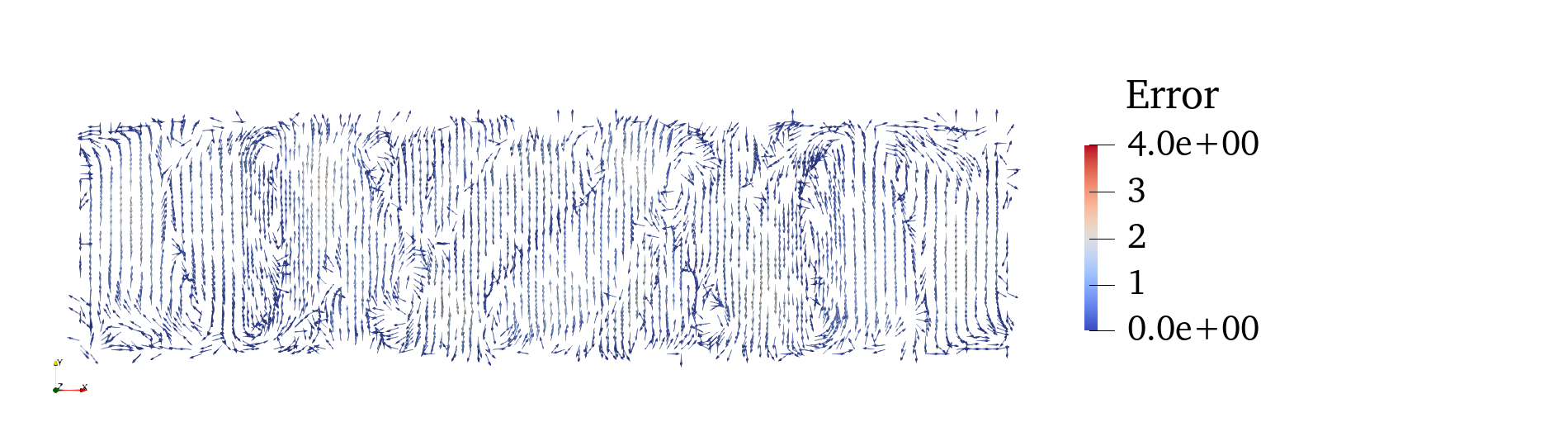}
\end{subfigure}
\begin{subfigure}{.5\textwidth}
    \centering
    \includegraphics[width=\textwidth]{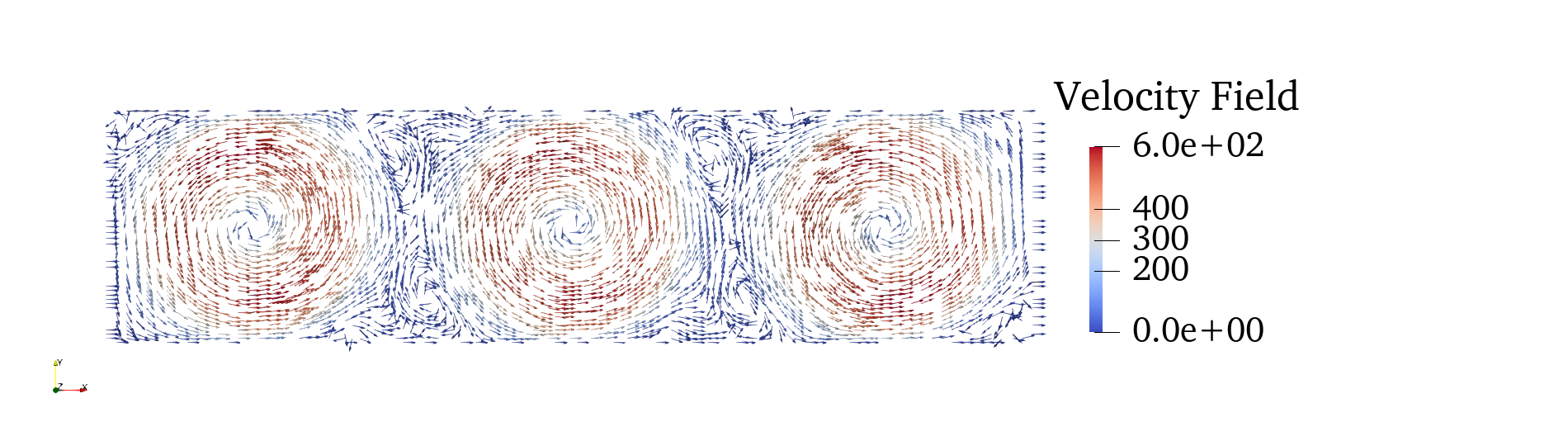}
\end{subfigure}
\begin{subfigure}{.5\textwidth}
    \centering
    \includegraphics[width=\textwidth]{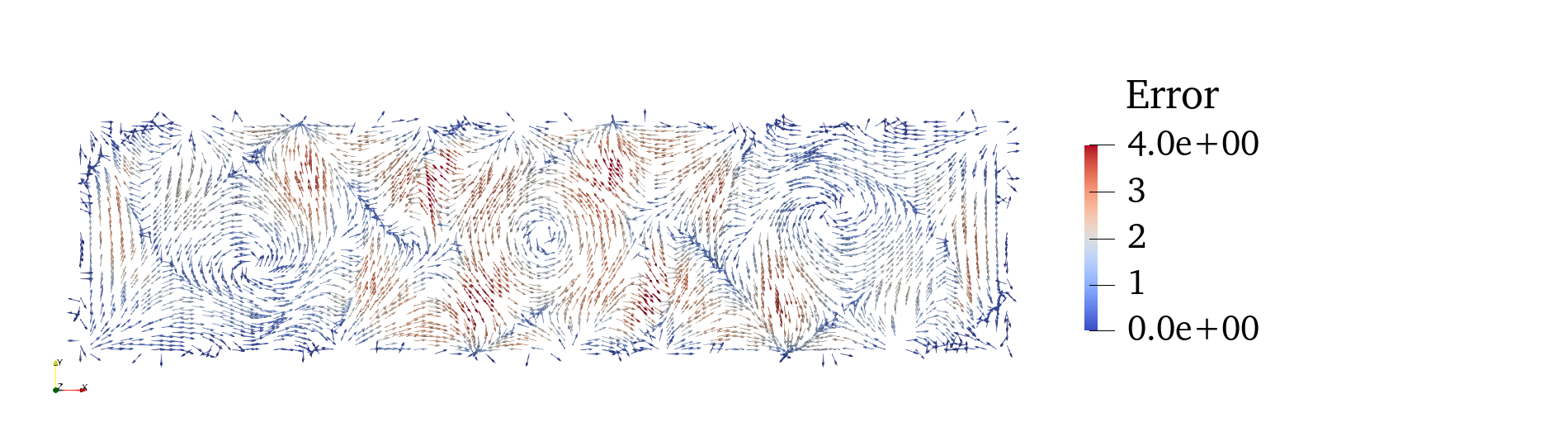}
\end{subfigure}
\caption{\footnotesize On the left, the approximated velocity field computed with our reduced order model for $Gr=132755$. On the right, the error committed with respect to the full order model solution. Both the columns are computed at different time-steps, respectively at $t=20000, t=50000, t=100000, t=200000$}
\label{fig:paraview}
\end{figure}

Finally, we report a visual example of some modes' coefficients predictions compared with the exact evolution (see Figure \ref{fig:modes_comp}) in the case of $Gr=132.755\cdot 10^3$, and the correspondent velocity field evolution in time (with the correspondent error committed) in Figure \ref{fig:paraview}. A frequency analysis is also visualized in Figure \ref{fig:FFT_modes_comp}, where the Fourier Transform of some modes coefficients evolution in time is reported. As it is noticeable in both Figure \ref{fig:modes_comp} and Figure \ref{fig:FFT_modes_comp} by the overlap of the red dashed-lines (exact coefficients and Fourier Transforms) and the coloured solid ones (approximated ones), the architecture succeeds in general in the prediction of all the coefficients dynamics. 

\newpage

\section{Conclusions and further developments}

In this work, we presented a novel approach to parametric time-dependent problems, consisting in a previous k-means clustering of the available solutions with respect to their associated parameters value. In this way, we proceeded with the training of $k$ C-LSTM independent models with the purpose of obtaining $k$ local representations of the solutions space. Each of these models was in principle able to generalize the problem's solution in a neighborhood of the parameters value with whom it had been trained, thus the second stage of the architecture was designed to find a non-linear function that combined in a proper way the predictions coming from the first-stage models.

This C-LSTM architecture has firstly been tested on low-dimensional ODEs systems such as the Duffing Oscillator parametrized in its non-linear component, and the Predator-Prey system. 

Subsequently, we presented promising results obtained in the case of the Rayleigh-Bérnard cavity flow, where the Incompressible Navier-Stokes equations in a rectangular cavity were considered. Here, a previous Proper Orthogonal Decomposition was applied to the discretized system, in order to properly reduce the problem's dimensionality. 

The results obtained were extremely positive considering the limited error propagation and the time reduction. Indeed, referring to the numerical solver taken as a baseline, the time needed for the online phase was decreased of the $97\%$.

\begin{center}
\begin{tabular}{ |c||c c| } 
 \hline
 \textbf{Netkar solver} & \textbf{Two stages architecture} & \\
 \hline
 \emph{Solving time} & \emph{Training time} & \emph{Online phase} \\ 
 \hline
 $\approx 3h$ & $\approx 4h$ & $\approx 5m$ \\ 
 \hline
\end{tabular}
\end{center}

It can finally be concluded that this method could reveal itself particularly useful in parametric large-scale systems, whose dynamic exhibits a non-linear behaviour difficult to generalize in the whole parameter space. 

Therefore, further applications are to be investigated for instance in the field of bifurcating systems, where such \emph{partitioning-averaging} approach could reveal crucial for a better and faster evaluation of the solution qualitative behaviour depending on the interplay between a given set of parameters.

It is in addition to be noted that the possibility of connecting multiple architectures with each other offers various advantages when it comes to different behaviours in time, also related to projection-based ROMs. As in the case of Section \ref{sec:res} indeed, where different PODs were performed in the initial \emph{swing-in} phase and in the \emph{periodic} one, different projection-based ROMs could be used to create multiple reduced-order spaces, later connected together with the evolution predicted by subsequent pre-trained two-stages architectures (one in every reducted order space).
\section*{Statements and Declarations}
\begin{itemize}
    \item \emph{Conflict of interest/Competing interests}
    
    The authors have no conflicts of interest or competing interests.
\end{itemize}
\newpage

\section{Appendix-A}\label{appendix}
\subsection{Recurrent Neural Networks and Long-Short Term Memory cells}

\begin{figure}[h!]
\centering
\input{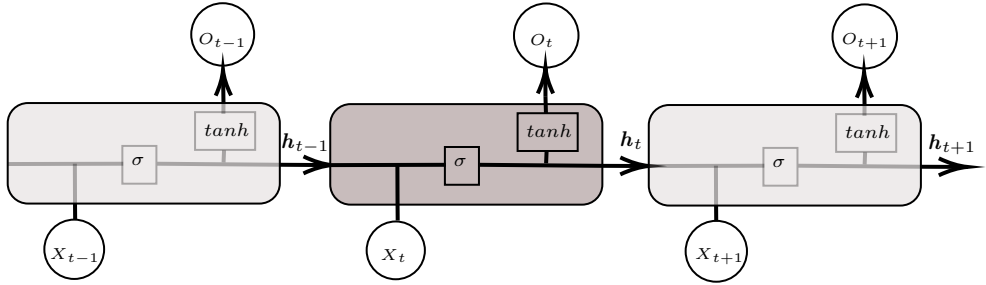}
\caption{\small The first historical example of a recurrent cell (ERNN).}
\label{fig:elman}
\end{figure}

Recurrent Neural Networks (RNN) are currently the most commonly used Neural Network architecture for sequence prediction problems \cite{RNN_review}. Every RNN is a combination of a certain number of RNN cells, which can be chosen among different realizations varying in complexity. However, all of them still carry out the same basic idea displayed in Figure \ref{fig:elman} and initially introduced by Elman \cite{ELMAN} in 1990.

He essentially proposes to implement a system of internal \emph{gates} aimed to build a bridge between the input and the output. In particular, this relation is mediated by a \emph{hidden state} (\emph{context cell}) $h_t$, managed through the equation (\ref{eq:2}) at each time-step, according to a trainable combination of the current input $x_t$ and the previous hidden state $h_{t-1}$. On the other hand, the cell's output $o_t$ at each temporal step is obtained through the equation (\ref{eq:3}), working on the current $h_t$.
\begin{align}
    \label{eq:2}
    h_t =& \sigma(W_h\cdot[h_{t-1},x_t]+b_h) \\
    \label{eq:3}
    o_t =& tanh(W_o\cdot h_{t}+b_o)
\end{align}
Here, the weight matrices $W_h$ and $W_o$, and the bias vectors $b_h$ and $b_o$ represent the trainable parameters of the network. 

It can be easily seen that those kind of update-laws insert feedback loops in the RNN cell, connecting its current state to the next one. These connections are of extreme importance in order to consider past information when updating the current cell state, conferring to the Recurrent Neural Network the possibility to preserve a \emph{memory of the system}.

However, the Elman Recurrent Neural cell suffers from the \emph{vanishing gradient} and \emph{exploding gradient} problems over very long sequences. This implies that the simple RNN cells are not capable of carrying long-term dependencies to the future: the back-propagated gradients tend to \emph{vanish} (and consequently the weights are not updated adequately) \cite{vanishing}, or \emph{explode} (resulting in unstable weight matrices). 

Over the years lots of variations have been proposed to overcome such problems. One of the currently more popular solutions are the Long-Short-Term-Memory cells (LSTM), first introduced in 1997 \cite{LSTM_intro}.

They present a more complex internal structure: 

\begin{align}
    q_t =& tanh(W_q\cdot[h_{t-1},x_t]+b_q) \\
    i_t =& \sigma(W_i\cdot[h_{t-1},x_t]+b_i) \\
    f_t =& \sigma(W_f\cdot[h_{t-1},x_t]+b_f) \\
    o_t =& \sigma(W_o\cdot[h_{t-1},x_t]+b_o) \\
    c_t =& f_t\odot c_{t-1} + i_t\odot q_t \\
    h_t =& o_t\odot tanh(c_t)
\end{align}

where $W_i$ , $b_i$ , $W_f$ , $b_f$ , $W_o$ , $b_o$ , $W_c$ , $b_c$  are the trainable weight matrices and bias vectors, while $\odot$ is the Hadamard product. 

The \emph{vanishing} and \emph{exploding gradient} problems have been solved with the introduction of Constant Error Carousels units (CECs) \cite{CECs}. Indeed, they enforce in the LSTM cells a system of internal gates and loops that makes them able to learn time lags of more than $1000$ discrete time steps, in contrast to previous ERNNs, which were already failing with time lags of $10$ time steps \cite{DL_NN}. From here, the name itself of LSTM cells is derived, underlying that they are able to capture both the short and the long-term dependencies in the training inputs.

\begin{figure}[t]
    \centering
    \input{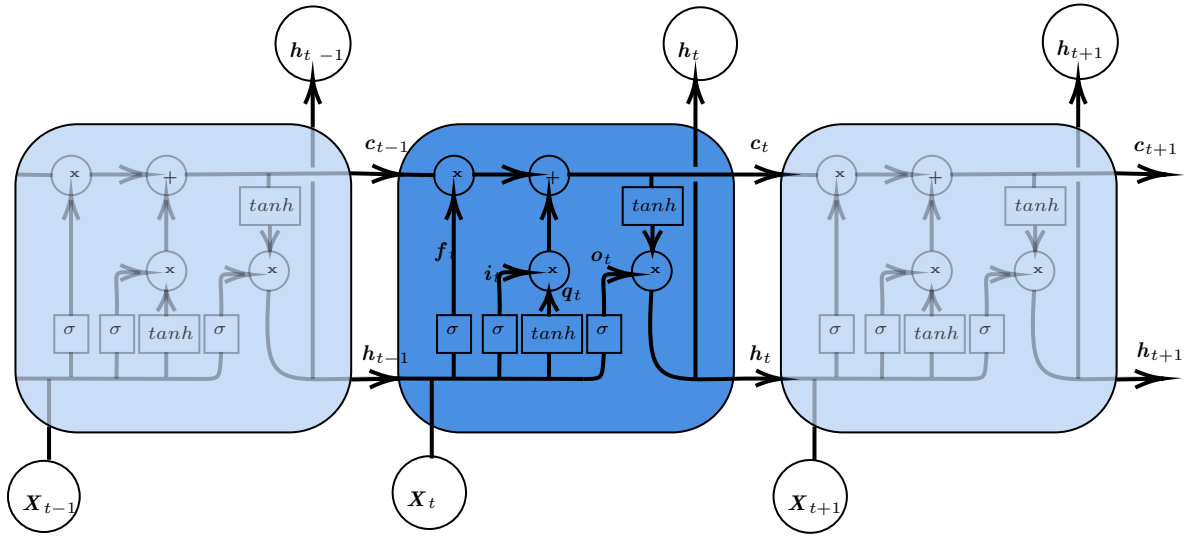}
    \caption{\small Internal logic structure of LSTM cells.}
    \label{fig:LSTM}
\end{figure}

Intuitively, such cells retain information about their past history through two quantities. Firstly $c_t$, which can be seen as the \emph{long-term memory} of the cell, and whom update is both designed for forgetting something of the past and incorporating new information coming from the current input. In second place, there is $h_t$, the \emph{hidden state} (and the output itself) of the cell, representing the \emph{short-term memory} component, updated both considering a non-linear transformation of the \emph{long-term memory} information $c_t$, and the \emph{output gate} one $o_t$.

The other quantities computed inside the LSTM cells can be explained as an interplay of gated structures, which are trained combining in a non-linear way the cell's \emph{hidden state} and \emph{input} information. Indeed, we can identify the \emph{input gate} $i_t$, aimed to perform a non-linear transformation of the current input, the \emph{forget gate} $f_t$ holding indications on the amount of past information which is safe to forget, and the \emph{output gate}, that offers a first proposal about the cell's final output.

\newpage
\section{Acknowledgements}
This work was partially funded by European Union Funding for Research and
Innovation --- Horizon 2020 Program --- in the framework of European Research
Council Executive Agency: H2020 ERC CoG 2015 AROMA-CFD project 681447 ``Advanced
Reduced Order Methods with Applications in Computational Fluid Dynamics'' P.I.
Professor Gianluigi Rozza. We also acknowledge the PRIN 2017 “Numerical Analysis
for Full and Reduced Order Methods for the efficient and accurate solution of
complex systems governed by Partial Differential Equations” (NA-FROM-PDEs).
\bibliographystyle{abbrv}
\bibliography{bib/biblio}
\end{document}